\documentclass[12pt]{article}     

\begin{document}
\begin{flushright} 
{ q-alg/0401088}
\end{flushright}

\begin{center}
 {\bf \huge  Cayley--Klein  Contractions  
of Quantum Orthogonal Groups in Cartesian Basis }         
\end{center}

\begin{center}
 {\bf N.A. Gromov, V.V. Kuratov }\\
 Department of Mathematics, 
Syktyvkar Branch of IMM UrD RAS,  \\
Chernova st., 3a, Syktyvkar, 167982, Russia \\
E-mail: gromov@dm.komisc.ru
\end{center}

\bigskip

\begin{center}
 {Abstract}
\end{center}
\noindent
Spaces of constant curvature and their motion groups
are described most naturally in Cartesian basis. 
All these motion groups also known as CK groups are obtained from
orthogonal group by contractions and analytical continuations.
On the other hand quantum deformation of orthogonal
group $SO(N) $ is most easily performed in so-called
symplectic basis. We reformulate its standard
quantum deformation  to Cartesian basis and
obtain all possible contractions of quantum orthogonal  group $SO_q(N)$
both for untouched and transformed deformation parameter.
It  turned out, that similar to  undeformed case all CK contractions
of $SO_q(N)$ are realized. An algorithm  for obtaining  nonequivalent 
(as Hopf algebra) contracted quantum groups is suggested.
Contractions of $SO_q(N), \; N=3,4,5$ are regarded as an examples.

\section{Introduction}     

Systematic definitions of quantum deformations of classical simple
Lie groups and algebras as well as  descriptions their properties
was given in  \cite{Fad-89}. 
Simple Lie groups and algebras are transformed by the
 contraction operation first introduced by E. Wigner and E. In{\"o}n{\"u} 
\cite{IW-53}  to a nonsemisimple ones.
Quantum analogues of the nonsemisimple low dimensions  Lie algebras 
was obtained by contractions of quantum algebras  $so_q(3), su_q(2), $
\cite{VK}--\cite{V-93} and contractions of low dimensions  quantum groups
was discussed in \cite{C2}--\cite{Gr-93}.
 Two types of contractions was discovered: with untouched deformation
parameter (in \cite{VK},\cite{SWZ} for quantum algebras  and
in \cite{G-93},\cite{Gr-93} for  quantum  groups)  and with transformed
deformation  parameter  \cite{C0},\cite{C1},\cite{V-93},\cite{C2}. 
For the last case the quantum  deformations of the algebras of 
the maximal symmetric motion groups of the $N-$dimensional flat spaces 
was  constructed in  \cite{V-95}.
$\gamma-$Poincare quantum group  was obtained by contractions of 
the orthogonal quantum group $ SO_q(N) $ \cite{Z}.  
Quantum Euclid group $E_{\kappa}(2)$ was described 
both by  contraction of $SU_q(2)$ \cite{S} and by direct quantization 
of Lie--Poisson structure \cite{M}. 
A separate line of investigation is presented by the $R$-matrix approach
to the quantum analogues of Euclid, Heisenberg and inhomogeneous
groups \cite{SWW}--\cite{AC}.

It is well known \cite{G-90} that the motion groups of all $ 3^{N-1} $
$(N-1)$-dimensional constant curvature spaces may be obtained
by contractions and analytic continuations of the classical orthogonal 
group $SO(N).$ Cayley--Klein groups is the short name for this set of groups.
The fundamental orthogonal $A^tA=I$ matrix $A \in SO(N)$
is replaced by the matrix $A(j)$ whose  elements
$(A(j))_{kp}=(k,p)a_{kp}, \;
(k,p)=\prod_{l=min\{k,p\}}^{max\{k,p\}-1}j_l, \; k,p=1,\dots,N  $
are subject of the additional $j$-orthogonality relations
$A(j)^tA(j)=1,$ where the parameters $j_k$ takes three values each
 $j_k=1,\iota_k,i.$ The commutative  
$\iota_k\iota_p=\iota_p\iota_k\neq 0,\; k\neq p $
nilpotent  $ \iota_k^2=0 $ units $ \iota_k $ are corresponded to contractions
and the imaginary unit $i^2=-1$ to analytic continuations.

In the case of the quantum orthogonal group $SO_q(N)$ additionally 
the deformation parameter $q=\exp{z}$ is transformed as follows \cite{S-97}:
$z=Jv,\; J=(1,N), $ where $v$ is the new deformation parameter .
At the same time the quantum  group contractions with untransformed
deformation parameter are known \cite{G-93},\cite{Gr-93}.
For unification of both such cases  in one approach the concept
of different couplings of Cayley-Klein and Hopf structures was
suggested \cite{VG},\cite{T-95}. It is well known that quantum
groups are Hopf algebras and   Cayley-Klein  structure is
defined by the distribution of the contraction parameters $j$
among the elements of the generating matrix.
For the quantum orthogonal group in so-called "symplectic" basis
(where the invariant quadratic form for $q=1$ is defined by the
matrix $C_0$ with all null elements except units on the  secondary
diagonal) this concept was realized in \cite{S-00}--\cite{GKK}
by the substitution in standard machinery of quantum group 
the generating matrix
$T_{\sigma}(j)=D_{\sigma}A(j)D_{\sigma}^{-1}, \; D_{\sigma}=DV_{\sigma},$
where the matrix $D$ is the solution of the equation   $D^tC_0D=I$
and describe transformation from Cartesian basis to symplectic one.
The matrix $V_{\sigma},\; (V_{\sigma})_{ik}=\delta_{\sigma_i,k},$
where  $\sigma \in S(N)$ is a $N$ order  permutation, define the
distribution of the contraction parameters in $T_{\sigma}(j).$
In this case the transformation of the deformation parameter depend on
permutation  $ \sigma. $
All permutations which leads to untouched $(J=1)$ deformation parameter
and some permutations which correspond to transformed ones
are enumerated in \cite{S-00}--\cite{GKK}. The contracted
quantum groups $ SO_v(N;j;\sigma) $ in these papers were regarded as 
Hopf algebra over Pimenov algebra $D(\iota)$ generated by nilpotent 
commutative generators.
It turned out that not all  Cayley--Klein contractions
are admissible for quantum groups in this assumption which therefore
is too restrictive.

The main statement of the algebraic structures contraction method
is to take into account in all relations only principal parts
with respect to tending to zero contraction parameter and to neglect
all others. Therefore in this paper  in all relations of quantum group 
theory  only principal (complex) terms  are taken into account 
and all other terms with nilpotent multipliers are neglected.
Besides contractions of  orthogonal quantum groups $ SO_v(N;j;\sigma)$ 
are regarded in  more usual Cartesian basis. For untouched deformation
parameter results are the same as in \cite{S-00}--\cite{GKK} and for
all other permutations deformation parameter is multiplied by 
$ J=\displaystyle{\bigcup^n_{k=1}(\sigma_k,\sigma_{k'})},$ where $n$
is  integral part of $ N/2.$
The unification of multipliers
$(\sigma_k,\sigma_p)\cup (\sigma_m,\sigma_r) $
is understood as the first power product of all parameters $j_k$
which appear at least in one multiplier
$(\sigma_k,\sigma_p)$ or $(\sigma_m,\sigma_r).$ 
For example, $ (j_1j_2)\cup (j_2j_3)=j_1j_2j_3.$
It  turned out, that the full scheme of CK contractions
are realized for the quantum group $SO_q(N).$
Not all identically contracted quantum groups corresponding 
to different permutations $\sigma$  are nonisomorphic.
Quantum groups  isomorphism is connected with the notion of 
equivalent distributions of nilpotent parameters in generating matrix.
Nonisomorphic contracted  quantum groups are correspond in the first place 
to generating matrices with nonequivalent distributions of nilpotent
parameters and secondly to equivalent generating matrices but with
different transformations of deformation parameter $(J_1 \neq J_2).$
As an example quantum groups $ SO_v(3;j;\sigma)$ are considered in detail
and nonisomorphic contractions  are given for quantum groups
$ SO_v(N;j;\sigma), \; N=4,5.$ The russian version of this paper 
was published in \cite{S-03}.

\section{Definition of quantum group $ SO_v(N;j;\sigma)$}

Let us  start with an algebra $ {\bf D} \langle (U(j;\sigma))_{ik} \rangle $
of noncommutative polynomials of  $ N^2 $ variables, which are an elements of
generating matrix
$ (U(j;\sigma))_{ik}= (\sigma_i,\sigma_k)u_{\sigma_i\sigma_k }. $
Let us introduce  the transformation of the deformation parameter
$q=e^z $ as follows: $z=Jv,$ where $v$ is a new deformation parameter 
and  $ J $ is some product of  parameters $ j $ for the present unknown.
Let $ \tilde R_v(j), \tilde{C}_v(j) $ be matrices which are obtained from
$\tilde{R}_q, \tilde{C}$ respectively
by the replacement of deformation parameter $ z $ with $ Jv. $
The commutation relations of the generators
$ U(j;\sigma) $ are defined by
\begin{equation}
\tilde{R}_v(j)U_1(j;\sigma)U_2(j;\sigma) =
U_2(j;\sigma)U_1(j;\sigma)\tilde{R}_v(j),
\label{1*}
\end{equation}
where    
$$ 
U_1(j;\sigma)=U(j;\sigma) \otimes I, \;
U_2(j;\sigma)=I \otimes U(j;\sigma), \;
$$
$$
U(j;\sigma)= V_{\sigma}U(j)V_{\sigma}^{-1}, \;
(V_{\sigma})_{ik}=\delta_{\sigma_i k},
$$ 
$$ 
\tilde{R}_v(j) = (D \otimes D)^{-1}R_v(j)(D \otimes D), \quad
R_v(j)=R_q(z \rightarrow Jv), \quad
$$ 
$$ 
D^{-1}=\frac{1}{\sqrt{2}}
\left ( \begin{array}{ccc}
      I & 0 &   {\tilde C_0} \\
      0 & \sqrt{2} &  0 \\
      i{\tilde C_0} & 0  &  -iI
      \end{array} \right ),    \       N=2n+1,
$$ 
$ {\tilde C_0} $ is the $ n \times n $ matrix with 
all null elements except units on the  secondary
diagonal and the explicit form of the matrix $ \tilde{R}_q $
in Cartesian basis is given in Appendix 1.
The additional relations of $(v,j)$-orthogonality are hold
\begin{equation}
U(j;\sigma)\tilde{C}_v(j)U^t(j;\sigma)=\tilde{C}_v(j),\; \;
U^t(j;\sigma)\tilde{C}^{-1}_v(j)U(j;\sigma)=\tilde{C}^{-1}_v(j), \;\;
\label{2*}
\end{equation}
where   
$ C=C_0q^{\rho}, $ 
and     
$ \rho=diag(\rho_1, \ldots, \rho_N), \
(C_0)_{ik}=\delta_{i'k}, \  i,k=1, \ldots, N,\ i'=N+1-i, \ $
 that is  
$ (C)_{ik}=q^{\rho_{i'}}\delta_{i'k} $ and   
$ (C^{-1})_{ik}=q^{-\rho_i}\delta_{i'k}, $
$\tilde{C}_v(j)=D^{-1}C_v(j)(D^t)^{-1}, $
$$ 
(\rho_1, \ldots, \rho_N)=
\left \{ \begin{array}{ccc}
     (n-\frac{1}{2}, n-\frac{3}{2}, \ldots , \frac{1}{2},0,-\frac{1}{2},
     \ldots , -n+\frac{1}{2}), \; N=2n+1 \\
     (n-1, n-2, \ldots, 1,0,0,-1, \ldots, -n+1), \; N=2n.
     \end{array} \right.
$$ 

{\it  The quantum orthogonal Cayley--Klein group
$ SO_v(N;j;\sigma)$ 
is defined as the quotient algebra of
${\bf D} \langle (U(j;\sigma))_{ik} \rangle $
by relations (\ref{1*}),(\ref{2*}).}

Formally  $ SO_v(N;j;\sigma)$ is a Hopf algebra with the following
coproduct $\Delta,$ counit $\epsilon$ and antipode $ S: $
$$
\epsilon(U(j;\sigma))=I, \;\;
\Delta U(j;\sigma)=U(j;\sigma) \dot{\otimes} U(j;\sigma), \;\;
$$
\begin{equation}
S(U(j;\sigma))=\tilde{C}_v(j)U^t(j;\sigma) \tilde{C}^{-1}_v(j),
\label{3*}
\end{equation}
where $ (A \dot{\otimes} B)_{ik}=\sum_{p} A_{ip}\otimes B_{pk}.$
The explicit form of antipode is given in Appendix 2 and 
$(v;j)$-orthogonality in Appendix 3.

{\it Remark.} All relations for the quantum group $ SO_v(N;j;\sigma)$
may be obtained from the corresponding relations for $ SO_q(N)$ 
in Cartesian basis \cite{S-97} by replacement  $z \rightarrow Jv$
and $ u_{ik} \rightarrow (\sigma_i,\sigma_k)u_{\sigma_i\sigma_k }. $

\section{The basic theorem}

According to the algebraic structures contraction method
in all relations of the previous section for nilpotent values of $j$
only principal (complex) terms are taken into account and 
all other terms with nilpotent multipliers are neglected.
Relation is called {\it admissible}, if it is possible to  select 
a principal  terms. Otherwise relation is called {\it inadmissible}.
For example, equation $a+\iota_1b+\iota_2c = a_1+\iota_1d$
is admissible equation and is equivalent to $ a=a_1,$
whereas  equation
$\iota_1b+\iota_2c=\iota_1\iota_2d$
is inadmissible.

The formal definition of the quantum group
$ SO_{v}(N;j;\sigma) $ should be a real definition of contracted 
quantum group, if the proposed construction is a consistent 
Hopf algebra structure for  the principal  terms of all relations
under nilpotent values of some or all parameters $j$.
In other words, if all relations of the previous section are admissible.
The following theorem holds.

To prove the admisibility of the relations it is necessary to have their
explicit expressions. Such expressions are obtained for coproduct, counit,
antipode and $(v,j)$-orthogonality relations for arbitrary $N.$
Commutation relations (\ref{1*}) for generators of orthogonal quantum group
are written as an overdetermined equation system and its explicit solution
has obtained only for $N=3.$
The following theorem holds.

{\bf Theorem.} {\it
If the  commutation relations for generators are defined and
the deformation parameter is transformed as $z=Jv, \;
 J=\displaystyle{\bigcup^n_{k=1}(\sigma_k,\sigma_{k'})},$
then all Caley--Klein contractions of
quantum groups $SO_v(N;j;\sigma)  $  are allowed.}

{\it Proof.}
Let us prove consistency of our construction for most singular case
when all parameters $j$ are  nilpotent. Counit 
$\varepsilon(u_{\sigma_i\sigma_k})=0,\, i\neq k, \,
\varepsilon(u_{\sigma_k\sigma_k})=1,\, k=1,\ldots ,n $
do not restrict the values of $j.$
Multiplier
$C_{ikr}=(\sigma_i,\sigma_r)(\sigma_r,\sigma_k)
(\sigma_i,\sigma_k)^{-1} $ 
in coproduct
$\Delta(u_{\sigma_i\sigma_k})=
\displaystyle{\sum^N_{r=1}C_{ikr}u_{\sigma_i\sigma_r}\otimes
u_{\sigma_r\sigma_k}} $
is equal to 1, if $\sigma_i < \sigma_r <\sigma_k, $
is equal to $(\sigma_k,\sigma_r)^2,$ if
$\sigma_i < \sigma_k <\sigma_r, $
and is equal to $(\sigma_r,\sigma_k)^2,$ if
$\sigma_r < \sigma_k <\sigma_i, $   
therefore all expressions for coproduct are admissible
for nilpotent values of  all $j.$
Because of symmetry
$(\sigma_i,\sigma_k)=(\sigma_k,\sigma_i) $
it is sufficiently to examine the case  $\sigma_k <\sigma_i. $

Let us analyze antipode $S(U(j;\sigma))$ (see Appendix 2).
Terms
$(\sigma_k,\sigma_{k'})^{-1}{\rm \sinh }(Jv\rho_k),$ $ k=1,\ldots ,n,$
are appeared in expressions (\ref{a-6}) for $p=n+1-k. $
They are well defined if one take multiplier $J$ equal to
the first power product of all parameters $j_k$
which appear at least in one multiplier
$(\sigma_k,\sigma_{k'})^{-1}, \, k=1, \ldots ,n,$ that is
$J=\displaystyle{\bigcup^n_{k=1}(\sigma_k,\sigma_{k'})}.$
Let us verify that all expressions for antipode are admissible.
For nilpotent $J$ one have ${\rm \sinh } J=J, \, {\rm \cosh } J=1.$
Two types of multipliers are appeared in antipode:
$$
A_{kM}(\alpha)=J
\Bigg( \frac{(\sigma_k,\sigma_M)}{(\sigma_{k'},\sigma_M)}\Bigg)^\alpha , \quad
B_{kM}(\alpha)=J^2
\Bigg( \frac{(\sigma_k,\sigma_M)}{(\sigma_{k'},\sigma_{M'})}\Bigg)^\alpha ,
$$
where $k=1,\ldots ,n, \ M=1,\ldots , N, \ \alpha=\pm 1. $
All these multipliers are well defined for nilpotent values of $j.$
Since $\alpha=\pm 1,$ then without loss of generality one assume
$\sigma_k <\sigma_{k'}.$
For $A_{kM}(\alpha) $ there are three possibilities:
$ (i)   \ \sigma_k <\sigma_M <\sigma_{k'}, \
(ii)   \ \sigma_M \leq \sigma_k <\sigma_{k'}, \
(iii)  \ \sigma_k <\sigma_{k'} \leq \sigma_M. $
In the case $ (i)  \ A_{kM}(1)=(\sigma_k,\sigma_M)^2, \;
A_{kM}(-1)=(\sigma_M,\sigma_{k'})^2, $
in the case $ (ii)  \ A_{kM}(1)=1, \;
A_{kM}(-1)=(\sigma_k,\sigma_{k'})^2, $
in the case $ (iii),$ on the contrary, 
$ A_{kM}(1)=(\sigma_k,\sigma_{k'})^2, \;
A_{kM}(-1)=1. $
Multipliers $ B_{kM}(\alpha) $ all the more are well defined
in view of $J^2. $ In particular, for most unfavorable case
$\sigma_{M'} <\sigma_M <\sigma_k <\sigma_{k'} $ 
one have the fraction
$(\sigma_k,\sigma_M)(\sigma_{k'}\sigma_{M'})^{-1}=
(\sigma_M,\sigma_{M'})^{-1}(\sigma_k,\sigma_{k'})^{-1}, $
but $J^2 $ contain the multiplier
$(\sigma_M,\sigma_{M'})(\sigma_k,\sigma_{k'}), $
therefore $B_{kM}(1)$ remains nonsingular.
If $J=1, $ then  $A_{kM}(\alpha)=B_{kM}(\alpha)=1.$
Because of arbitrary choice of $k$ and $M $
multipliers $A_{kM}(\alpha) $ and $B_{kM}(\alpha) $
are well defined for all values of $k$ and $M. $

Besides Hopf structure the $(v,j)$-orthogonality relations
(\ref{o1-1}--\ref{o1-9}),(\ref{ort-1}) are imposed on generators
of quantum group $SO_v(N;j;\sigma). $ Equations
(\ref{o1-5}),(\ref{o1-1}) and (\ref{o1-2}) for $ k=p $  evidently are
admissible. Equations (\ref{o1-3}),(\ref{o1-4}) for $p=n+1-k $
after division of both parts on $ (\sigma_k,\sigma_{k'}) $ have
terms with multipliers $ C_{kMk'},$ which are equal to 1, if
$\sigma_k <\sigma_{M} <\sigma_{k'} $ and are some product of
$j_k^2$  otherwise. Therefore these equations are admissible 
The rest equations of $(v,j)$-orthogonality have terms with
coefficients
$$
A_{KPM}=
\frac{(\sigma_K,\sigma_M)(\sigma_M,\sigma_P)}
{(\sigma_K,\sigma_P)} , \quad
B_{KPr}=J
\frac{(\sigma_K,\sigma_r)(\sigma_P,\sigma_{r'})}
{(\sigma_K,\sigma_P)},
$$
where $K, P, M=1,\ldots , N, \, r=1,\ldots ,n. $
These coefficients are well defined for all nilpotent values
of $j.$  For $A_{KPM} $ it is easily follow from the analysis
of three possible cases:
$ (i) \ \sigma_K <\sigma_M <\sigma_P, \,
(ii) \ \sigma_M <\sigma_K <\sigma_P, \,
(iii) \ \sigma_K <\sigma_P <\sigma_M. $ 
Moreover in the case $(i) \ A_{KPM}=1 $ and corresponding terms are
complex. Nonsingularity  of $B_{KPr} $ follows from simple analysis
of three possible cases:
$(a) \ \sigma_K<\sigma_r<\sigma_P<\sigma_{r'}, \,
(b) \ \sigma_K<\sigma_P<\sigma_r<\sigma_{r'}, \,
(c) \ \sigma_K<\sigma_r<\sigma_{r'}<\sigma_P. $

Thus we conclude, that  $(v,j)$-orthogonality relations are admissible
for any permutations and for nilpotent values of any parameters,
therefore they do not restrict contractions of quantum group.

\section{Nonisomorphic contracted quantum groups}

If all parameters $j_k=1,$ then the map
$u_{ik} \rightarrow (\sigma_i,\sigma_k)u_{\sigma_i\sigma_k } $
is invertible and all quantum groups $ SO_v(N;j;\sigma)$ 
for any $\sigma \in S_N $ are isomorphic as Hopf algebras.
Nonisomorphic quantum groups may appear under contractions
when all or some parameters $j$ take nilpotent values.
It is clear that nonisomorphic quantum groups are appear under
contractions with different numbers of parameters.
Contractions on the same parameters, but with different
transformations of deformation parameter (with different $J$)
naturally give in result nonisomorphic quantum groups.
Isomorphic quantum groups may appear under contractions
of $ SO_v(N;j;\sigma)$ with different $\sigma $ by   equal
numbers of  parameters, when multiplier $J$ include  equal
numbers of  parameters (but not necessarily the same) or when $J=1.$ 
In our approach contractions of quantum groups (even on equal
numbers of  parameters) are distinguished by the distributions 
of nilpotent parameters $j$ in generating matrix $U({j;\sigma}).$
Really, all relations of quantum group theory (commutators,
$(v,j)$-orthogonality, antipode, coproduct and counit) depend on
permutation $\sigma$ by means of generating matrix, while
matrices $R_v(j), C_v(j)$ depend on  $\sigma$ via  transformations 
of deformation parameter, that is via $J.$ Isomorphism of
contracted quantum orthogonal groups is described by the 
following theorem.

{\bf Theorem.} {\it
Quantum  groups $ SO_v(N;j;\sigma_1)$ and $ SO_w(N;j;\sigma_2) $
are  isomorphic, if the following relations for their generators holds:
\begin{equation}
U({j;\sigma_1})=V_{\sigma} U({j;\sigma_2}) V^{-1}_{\sigma},
\label{preobr}
\end{equation}
where matrix $ V_{\sigma}, \; \sigma \in  S_N $ satisfy
\begin{equation}
(V_{\sigma}\otimes V_{\sigma}) \tilde{R}_w(j)
(V_{\sigma} \otimes V_{\sigma})^{-1} = \tilde{R}_v(j), \quad
V_{\sigma}\tilde{C}_w(j) V_{\sigma}^t= \tilde{C}_v(j)
\label{uslovie}
\end{equation}
for  $ w=\pm v$ and $J_1=J_2 $ 
with possible replacement
$j_k $ on $ j_{N-k}, \; k=1,\dots,N-1. $}

{\it Proof.}
Commutation relations (\ref{1*}) of  $ SO_v(N;j;\sigma_1) $
after transformation (\ref{preobr}) take the form
$$
\tilde{R}_v(j)(V_{\sigma}\otimes V_{\sigma}) U_1(j;\sigma_2)U_2(j;\sigma_2)
(V_{\sigma}\otimes V_{\sigma})^{-1}=
$$
$$
(V_{\sigma}\otimes V_{\sigma}) U_2(j;\sigma_2)U_1(j;\sigma_2)
(V_{\sigma}\otimes V_{\sigma})^{-1} \tilde{R}_v(j)
$$
or after left multiplying  on $(V_{\sigma}\otimes V_{\sigma})^{-1} $
and right multiplying  on $V_{\sigma}\otimes V_{\sigma}, $ in the form
$$
(V_{\sigma}\otimes V_{\sigma})^{-1} \tilde{R}_v(j) 
(V_{\sigma} \otimes V_{\sigma}) U_1(j;\sigma_2)U_2(j;\sigma_2)=
$$
$$
U_2(j;\sigma_2)U_1(j;\sigma_2)(V_{\sigma}\otimes V_{\sigma})^{-1} 
\tilde{R}_v(j) (V_{\sigma}\otimes V_{\sigma}),
$$
which give first equation in (\ref{uslovie}).
Antipode (\ref{3*}) after transformation (\ref{preobr}) take the form
$$
V_{\sigma}S(U(j;\sigma_2))V_{\sigma}^{-1}=
\tilde{C}_v(j)\left(V_{\sigma}^{-1}\right)^tU^t(j;\sigma_2)V_{\sigma}^t 
\tilde{C}^{-1}_v(j)
$$
or
$$
S(U(j;\sigma_2))=
V_{\sigma}^{-1}\tilde{C}_v(j)\left(V_{\sigma}^{-1}\right)^tU^t(j;\sigma_2)
V_{\sigma}^t \tilde{C}^{-1}_v(j)V_{\sigma}.
$$
The last equation is just antipode of $ SO_v(N;j;\sigma_2), $
if take into account the second equation in  (\ref{uslovie}).
 At last, $(v,j)$-orthogonality relations (\ref{2*}) after (\ref{preobr})
take the form
$$
V_{\sigma}U(j;\sigma_2)V_{\sigma}^{-1}\tilde{C}_v(j)\left(V_{\sigma}^{-1}\right)^tU^t(j;\sigma_2)V_{\sigma}^t=
\tilde{C}_v(j)
$$
or
$$
U(j;\sigma_2)V_{\sigma}^{-1}\tilde{C}_v(j)\left(V_{\sigma}^{-1}\right)^tU^t(j;\sigma_2) =
V_{\sigma}^{-1}\tilde{C}_v(j)\left(V_{\sigma}^t\right)^{-1},
$$
which evidently is condition  (\ref{uslovie}) for matrix  $\tilde{C}_v(j). $

As a consequence of theorem is the following algorithm of obtaining
of nonisomorphic contracted quantum groups.
One call two distributions of nilpotent parameters among elements
of generating matrices $U(j;\sigma_1), U(j;\sigma_2)$ {\it equivalent},
if they are connected by two operations:
1) they pass in each other by the permutations of the same columns and rows
   of generating matrices, that is by  (\ref{preobr});
2) matrices  pass in each other by reflection relative secondary diagonal
   with possible simultaneous replacement of $j_k $ with $ j_{N-k}, \; 
   k=1,\dots,N-1. $
Nonisomorphic contracted quantum groups are corresponded 
in the first place to the nonequivalent generating matrices
and secondly to equivalent generating matrices, but with
different transformations of deformation parameters $(J_1 \neq J_2).$
For illustration of algorithm all nonequivalent contractions 
of quantum groups $ SO_v(N;j;\sigma),\; N=3,4,5$ shall be regarded
in the next sections.

\section{Quantum groups $ SO_v(3;j;\sigma)$}

Quantum group $ SO_q(3)$ has four nonisomorphic contracted  groups:
two Euclid groups
$E_v^0(2)=SO_v(3;\iota_1,j_2;\sigma_0),\; J=\iota_1, $ 
$ E_z(2)=SO_z(3;\iota_1,1;\sigma), \; J=1,$
where $\sigma_0=(1,2,3),\; \sigma=(2,1,3),  $
and two Galilei groups
$ G_v^0(2)=SO_v(3;\iota_1,\iota_2;\sigma_0), \;J=\iota_1\iota_2, $ 
$G_v(2)=SO_v(3;\iota_1,\iota_2;\sigma),$ $J=\iota_2. $
For comparison, nondeformed complex rotation group $SO(3)$
has two nonisomorphic Cayley--Klein contracted  groups:
Euclid group $E(2)$ and  Galilei group  $G(2).$

\subsection{\bf Quantum groups
$ SO_v(3;j;\sigma_0),\; \sigma_0=(1,2,3) $}

Let $C_1={\rm \cosh } Jv,\; S_1={\rm \sinh } Jv,\; J=j_1j_2. $
Generating matrix
\begin{equation}
U(j)=\left (\begin{array}{ccc}
u_{11} & j_1u_{12}& j_1j_2u_{13} \\
j_1u_{21} & u_{22} & j_2u_{23} \\
j_1j_2u_{31} & j_2u_{32} & u_{33} \end{array} \right )
\label{a1}
\end{equation}
satisfy $(v,j)$-orthogonality relations:
(i) $ U(j)C_v(j)U^{t}(j)=C_v(j),$  i.e.
$$
iJS_1[u_{13},u_{11}]=C_1(u^2_{11}+J^2u^2_{13}-1)+j_1^2u_{12}^2,
$$
$$
iJS_1[u_{23},u_{21}]=C_1(j_1^2u^2_{21}+j^2_2u^2_{23})+u_{22}^2-1,
$$
$$
iJS_1[u_{33},u_{31}]=C_1(J^2u^2_{31}+u^2_{33}-1)+j_2^2u_{32}^2,
$$
$$
u_{11}u_{21}j_1C_1-iu_{13}u_{21}j_1JS_1+j_1u_{12}u_{22}+
u_{13}u_{23}j_2JC_1+iu_{11}u_{23}j_2S_1=0,
$$
$$
u_{11}u_{31}JC_1-iu_{13}u_{31}J^2S_1+Ju_{12}u_{32}+
u_{13}u_{33}JC_1+iu_{11}u_{33}S_1=iS_1,
$$
$$
u_{21}u_{31}j_1JC_1-iu_{23}u_{31}j_2JS_1+j_2u_{22}u_{32}+
j_2u_{23}u_{33}C_1+iu_{21}u_{33}j_1S_1=0,
$$
$$
u_{21}u_{11}j_1C_1-iu_{23}u_{11}j_2S_1+j_1u_{22}u_{12}+
u_{23}u_{13}j_2JC_1+iu_{21}u_{13}j_1JS_1=0,
$$
$$
u_{31}u_{11}JC_1-iu_{33}u_{11}S_1+Ju_{32}u_{12}+
u_{33}u_{13}JC_1+iu_{31}u_{13}J^2S_1=-iS_1,
$$
\begin{equation}
u_{31}u_{21}j_1JC_1-iu_{33}u_{21}j_1S_1+j_2u_{32}u_{22}+
u_{33}u_{23}j_2C_1+iu_{31}u_{23}j_2JS_1=0
\label{a2}
\end{equation}
and (ii) $U(j)^{t}C^{-1}_v(j)U(j)=C^{-1}_v(j),$ i.e.
$$
iJS_1[u_{11},u_{31}]=C_1(u^2_{11}+J^2u^2_{31}-1)+j_1^2u^2_{21},
$$
$$
iJS_1[u_{12},u_{32}]=C_1(j_1^2u^2_{12}+j^2_2 u^2_{32})+u^2_{22}-1,
$$
$$
iJS_1[u_{13},u_{33}]=C_1(u^2_{33}+J^2u^2_{13}-1)+j_2^2u^2_{23},
$$
$$
j_1u_{11}u_{12}C_1+iu_{31}u_{12}j_1JS_1+j_1u_{21}u_{22}+
Ju_{31}u_{33}C_1-iu_{11}u_{33}S_1=0,
$$
$$
Ju_{11}u_{13}C_1+iu_{31}u_{13}J^2S_1+Ju_{21}u_{23}+
Ju_{13}u_{33}C_1-iu_{11}u_{33}S_1=-iS_1,
$$
$$
j_1Ju_{12}u_{13}C_1+iu_{32}u_{13}j_2JS_1+j_2u_{22}u_{23}+
j_2u_{32}u_{33}C_1-iu_{12}u_{33}j_1S_1=0,
$$
$$
j_1u_{12}u_{11}C_1+iu_{32}u_{11}j_2S_1+j_1u_{22}u_{21}+
j_2Ju_{32}u_{31}C_1-iu_{12}u_{31}j_1JS_1=0,
$$
$$
Ju_{13}u_{11}C_1+iu_{33}u_{11}S_1+Ju_{23}u_{21}+
Ju_{33}u_{31}C_1-iu_{13}u_{31}J^2S_1=iS_1,
$$
\begin{equation}
j_1Ju_{13}u_{12}C_1+iu_{33}u_{12}j_1S_1+j_2u_{23}u_{22}+
j_2u_{33}u_{32}C_1-iu_{13}u_{32}j_2JS_1=0.
\label{a3}
\end{equation}
There are three independent generators, for example,
$ u_{12}, u_{13}, u_{23},$  which are situated above diagonal.
Their commutators are obtained from $RUU$-relations
$\tilde{R}_v(j)U_1(j)U_2(j)=U_2(j)U_1(j)\tilde{R}_v(j)$
and are in the form 
$$
[u_{12},u_{23}]=i\frac{{\rm \sinh } Jv}{J}u_{22}(u_{11}-u_{33}),
$$
$$
[u_{13},u_{23}]=u_{23}\left\{({\rm \cosh } Jv-1)u_{13}-
i\frac{{\rm \sinh } Jv}{J} u_{33}\right\},
$$
\begin{equation}
[u_{12},u_{13}]=\left\{({\rm \cosh } Jv-1)u_{13}+
i\frac{{\rm \sinh } Jv}{J} u_{11}\right\}u_{12}.
\label{a4}
\end{equation}

An associative algebra $ SO_v(3;j;\sigma_0) $ is Hopf algebra
with counit  $ \epsilon(U(j))=I,$ i.e.
$\epsilon(u_{ik})=0, \epsilon(u_{kk})=1, $
coproduct  $\Delta U(j)=U(j)\otimes U(j)  $ in the form
$$
\Delta u_{12}=u_{11}\otimes u_{12}+
u_{12}\otimes u_{22}+j^2_2u_{13}\otimes u_{32}, \;
\Delta u_{21}=u_{21}\otimes u_{11}+
u_{22}\otimes u_{21}+j^2_2u_{23}\otimes u_{31},
$$
$$
\Delta u_{23}=u_{22}\otimes u_{23}+
u_{23}\otimes u_{33}+j_1^2u_{21}\otimes u_{13},\;
\Delta u_{32}=u_{32}\otimes u_{22}+
u_{33}\otimes u_{32}+j_1^2u_{31}\otimes u_{12},
$$
$$
\Delta u_{13}=u_{11}\otimes u_{13}+
u_{12}\otimes u_{23}+u_{13}\otimes u_{33}, \;
\Delta u_{31}=u_{31}\otimes u_{11}+
u_{32}\otimes u_{21}+u_{33}\otimes u_{31},
$$
$$
\Delta u_{11}=u_{11}\otimes u_{11}+
j_1^2u_{12}\otimes u_{21}+J^2u_{13}\otimes u_{31},\;
$$
$$
\Delta u_{22}=u_{22}\otimes u_{22}+
j_1^2u_{21}\otimes u_{12}+j^2_2u_{23}\otimes u_{32},
$$
\begin{equation}
\Delta u_{33}=u_{33}\otimes u_{33}+
j_2^2u_{32}\otimes u_{23}+J^2u_{31}\otimes u_{13},
\label{a5}
\end{equation}
and antipode 
$ S(u(j))=C_v(j)U^t(j)C^{-1}_v(j),  $   where
$$
S(u_{12})=u_{21}{\rm \cosh }(\frac{Jv}{2})+ij_2^2u_{23}
\frac{1}{J}{\rm \sinh }(\frac{Jv}{2}), \;
$$
$$
S(u_{21})=u_{12}{\rm \cosh }(\frac{Jv}{2})+ij_2^2u_{32}
\frac{1}{J}{\rm \sinh }(\frac{Jv}{2}),
$$
$$
S(u_{23})=u_{32}{\rm \cosh }(\frac{Jv}{2})-ij_1^2u_{12}\frac{1}{J}
{\rm \sinh }(\frac{Jv}{2}),\;
$$
$$
S(u_{32})=u_{23}{\rm \cosh }(\frac{Jv}{2})-ij_1^2u_{21}\frac{1}{J}
{\rm \sinh }(\frac{Jv}{2}),
$$
$$
S(u_{13})=u_{31}{\rm \cosh }^2(\frac{Jv}{2})+u_{13}{\rm \sinh }^2(\frac{Jv}{2})+
i\frac{1}{2}(u_{33}-u_{11})\frac{1}{J}{\rm \sinh }(Jv),
$$
$$
S(u_{31})=u_{13}{\rm \cosh }^2(\frac{Jv}{2})+u_{31}{\rm \sinh }^2(\frac{Jv}{2})+
i\frac{1}{2}(u_{33}-u_{11})\frac{1}{J}{\rm \sinh }(Jv),
$$
$$
S(u_{11})=u_{11}{\rm \cosh }^2(\frac{Jv}{2}) -u_{33}{\rm \sinh }^2(\frac{Jv}{2})
+i\frac{1}{2}(u_{13}+u_{31})J{\rm \sinh }(Jv),
$$
$$
S(u_{33})=u_{33}{\rm \cosh }^2(\frac{Jv}{2}) -u_{11}{\rm \sinh }^2(\frac{Jv}{2})
-i\frac{1}{2}(u_{13}+u_{31})J{\rm \sinh }(Jv),
$$
\begin{equation}
\quad S(u_{22})=u_{22}.
\label{a6}
\end{equation}

{\it Remark.} Coproduct and counit of $ SO_v(3;j;\sigma) $
are the same for any permutation   $\sigma. $ Only antipode,
commutation and $(v;j)$-orthogonality relations are depend on $\sigma. $

For $ j_1=\iota_1 $   {\bf quantum Euclid group 
$E_v^0(2)=SO_v(3;\iota_1,j_2;\sigma_0), \; J=\iota_1 $ } is obtained.
From $(v;j)$-orthogonality relations it follows
$u_{11}=1,\; u_{22}=u_{33}, \; u_{23}=-u_{32}, $
and from   $RUU$-equations it follows that all these generators commute
and generate rotation group $SO(2).$ Therefore it is naturally 
to introduce new notations
$ u_{22}=u_{33}=\cos \varphi, \; u_{23}=\sin \varphi=-u_{32}, $
and rewrite generating matrix as
\begin{equation}
U(\iota_1;\sigma_0)=\left (\begin{array}{ccc}
1 & \iota_1u_{12} & \iota_1u_{13} \\
\iota_1u_{21} & \cos\varphi & \sin\varphi \\
\iota_1u_{31} & -\sin\varphi & \cos\varphi \end{array} \right )
\sim \left (\begin{array}{ccc}
\cdot & \circ & \circ \\
      & \cdot &  \cdot \\
 &  & \cdot \end{array} \right ),
\label{a7}
\end{equation}
where from $(v;j)$-orthogonality relations it follows
$$
u_{21}=-(u_{12}\cos\varphi+u_{13}\sin\varphi+i\frac{v}{2}\sin\varphi),\;\;
$$
\begin{equation}
u_{31}=u_{12}\sin\varphi-u_{13}\cos\varphi+i\frac{v}{2}(1-\cos\varphi).
\label{a8}
\end{equation}
Here and later
distribution of nilpotent parameters among elements of generating matrix
is shown with the help of notations:
$\circ = \iota_1,\;\bullet =\iota_2, \;\times = \iota_1\iota_2. $
(Let us remind that this distribution is symmetric relatively diagonal).
Dots denote complex elements.
Commutation relations of independent generators are as follows
$$
[u_{12},\sin \varphi ]=iv\cos \varphi(1-\cos \varphi),\quad
$$
\begin{equation}
[\sin \varphi,u_{13}]=iv\sin \varphi\cos \varphi,\quad
[u_{12},u_{13}]=ivu_{12}.
\label{a9}
\end{equation}
Coproduct of quantum Euclid group is given by
$$
\Delta u_{12}=
1\otimes u_{12}+
u_{12}\otimes \cos\varphi-j_2^2u_{13}\otimes \sin\varphi,
$$
$$
\Delta u_{13}=1\otimes u_{13}+u_{12}\otimes \sin\varphi+
u_{13}\otimes\cos\varphi,
$$
\begin{equation}
\Delta\sin\varphi=
\cos\varphi\otimes\sin\varphi+\sin\varphi\otimes\cos\varphi,\;\;
\Delta\varphi=1\otimes\varphi+\varphi\otimes1,
\label{a10}
\end{equation}
antipode is as follows
$$
S(u_{12})=
-u_{12}\cos\varphi-u_{13}\sin\varphi, \quad
S(u_{13})=
-u_{13}\cos\varphi+u_{12}\sin\varphi, \;
$$
\begin{equation}
S(\varphi)=-\varphi,
\label{a11}
\end{equation}
and their counit is equal to zero:
$\epsilon (u_{12})=\epsilon (\varphi)=\epsilon (u_{13})=0.$

If   $u_{21},\;u_{31},\;\varphi$ are taken as independent generators,
then equations (\ref{a8})--(\ref{a11}) are rewritten
in the following way:
from $(v;j)$-orthogonality relations
$$
u_{12}=-u_{21}\cos\varphi+u_{31}\sin\varphi-i\frac{v}{2}\sin\varphi,
$$
\begin{equation}
u_{13}=-u_{21}\sin\varphi-u_{31}\cos\varphi-i\frac{v}{2}(1-\cos\varphi),
\label{a8.1}
\end{equation}
commutation relations
$$
[u_{21},\sin \varphi ]=iv\cos \varphi(1-\cos \varphi),\quad
$$
\begin{equation}
[\sin \varphi,u_{31}]=-iv\sin \varphi\cos \varphi,\quad
[u_{31},u_{21}]= ivu_{21},
\label{a9.1}
\end{equation}
coproduct
$$
\Delta u_{21}=
u_{21}\otimes 1+
\cos\varphi\otimes u_{21}+\sin\varphi\otimes u_{31},
$$
$$
\Delta u_{31}=u_{31}\otimes 1-\sin\varphi \otimes u_{21}+
\cos\varphi \otimes u_{31},
$$
\begin{equation}
\Delta\varphi=1\otimes\varphi+\varphi\otimes1,
\label{a10.1}
\end{equation}
antipode
$$
S(u_{21})=
-u_{21}\cos\varphi+u_{31}\sin\varphi -iv\sin \varphi, \quad
$$
\begin{equation}
S(u_{31})=
-u_{31}\cos\varphi-u_{21}\sin\varphi+iv(\cos \varphi -1), \quad
S(\varphi)=-\varphi
\label{a,b}
\end{equation}
and counit
$\epsilon (u_{21})=\epsilon (\varphi)=\epsilon 
(u_{31})=0.$

 Under contraction $ j_2=\iota_2  $  {\bf quantum analog 
$N_v^0(2)=SO_v(3;j_1,\iota_2;\sigma_0), \; J=\iota_2 $
of cylindrical group or Newton group} $N(2)$ is obtained.
Similarly to previous case with the help of $(v,j)$-orthogonality
relations the generating matrix may be written in the form
\begin{equation}
U(\iota_2;\sigma_0)=\left (\begin{array}{ccc}
\cos\psi & \sin\psi & \iota_2u_{13} \\
-\sin\psi & \cos\psi & \iota_2u_{23} \\
\iota_2u_{31} & \iota_2u_{32} & 1 \end{array} \right )
\sim\left (\begin{array}{ccc}
\cdot & \cdot & \bullet \\
 & \cdot & \bullet \\
 &  & \cdot \end{array} \right ),
\label{a12}
\end{equation}
where
$$
u_{31}= u_{23}\sin\psi- u_{13}\cos\psi +i{v \over 2}
(1-\cos\psi),\;
$$
\begin{equation}
u_{32}=- u_{23}\cos\psi -  u_{13}\sin\psi
-i{v\over 2} \sin\psi,
\label{a13}
\end{equation}
and independent generators are subject of commutation relations
$$
[\sin\psi,u_{23}]=iv\cos\psi(\cos\psi-1), \;
$$
\begin{equation}
[u_{23},u_{13}]=ivu_{23},\quad
[\sin\psi,u_{13}]=iv\sin\psi \cos\psi.
\label{a14}
\end{equation}
Hopf algebra is defined by coproduct
$$
\Delta(\sin\psi)=
\cos\psi\otimes\sin\psi+\sin\psi\otimes\cos\psi,\;\;
\Delta(\psi)=1\otimes\psi+\psi\otimes1,
$$
$$
\Delta u_{13}=u_{13}\otimes 1+\cos\psi\otimes u_{13}+
\sin\psi\otimes u_{23},
$$
\begin{equation}
\Delta u_{23}=
u_{23}\otimes 1+\cos\psi\otimes u_{23}-j_1^2\sin\psi\otimes u_{13},
\label{a15}
\end{equation}
by antipode
$$
S(u_{13})=
u_{31}+i\frac{v}{2}(u_{33}-u_{11})=
u_{23}\sin\psi-u_{13}\cos\psi+iv(1-\cos\psi),
$$
\begin{equation}
S(u_{23})=
u_{32}-i\frac{v}{2}j_1^2u_{12}=
-u_{23}\cos\psi-u_{13}\sin\psi-iv\sin\psi, \quad
S(\psi)=-\psi,
\label{a16}
\end{equation}
and by counit
$\epsilon (\psi)=\epsilon (u_{13})=\epsilon (u_{23})=0.$

The distribution of $\iota_1$ in matrix (\ref{a7}) is passed to
the distribution of $\iota_2$ in matrix (\ref{a12}) under reflection
on secondary diagonal and simultaneous substitution  $ J=\iota_1$
by $ J=\iota_2. $  This means that the quantum Euclid group
$E_v^0(2)=SO_v(3;\iota_1,1;\sigma_0) $  is isomorphic to
the quantum Newton group $N_v^0(2)=SO_v(3;1,\iota_2;\sigma_0) $
as well as in nondeformed case.
Under substitution
$u_{31}$ on $u_{13},$ $u_{21}$ on  $u_{23},$ 
$\varphi$ on $-\psi,$ $v$ on $-v$
commutation relations  (\ref{a9.1}) are transformed in (\ref{a14}),
coproduct (\ref{a10.1}) is transformed in  (\ref{a15}) and antipode
 --- in (\ref{a16}).

Two-dimensional contraction
$ j_1=\iota_1, \, j_2=\iota_2 $  gives {\bf quantum Galilei group}
$ G_v^0(2)=SO_v(3;\iota_1,\iota_2;\sigma_0), \; J=\iota_1\iota_2. $
With the help of $(v;j)$-orthogonality relations the generating matrix
may be written in the form
\begin{equation}
U(\iota;\sigma_0)=\left (\begin{array}{ccc}
1 & \iota_1u_{12} & \iota_1\iota_2u_{13} \\
-\iota_1u_{12} & 1 & \iota_2u_{23} \\
\iota_1\iota_2u_{31} & -\iota_2u_{23} & 1 \end{array} \right )
\sim\left (\begin{array}{ccc}
\cdot & \circ & \times \\
 & \cdot & \bullet \\
 &  & \cdot \end{array} \right ),
\label{a17}
\end{equation}
where
$ u_{31}=-u_{13}+u_{12}u_{23}, $
and independent generators satisfy commutation relations
\begin{equation}
[u_{12},u_{23}]=0, \quad
[u_{23},u_{13}]=ivu_{23}, \quad
[u_{12},u_{13}]=ivu_{12}.
\label{a18}
\end{equation}
Hopf algebra structure is given by coproduct
$$
\Delta u_{12}=
1\otimes u_{12}+u_{12}\otimes 1, \quad
\Delta u_{23}=
1\otimes u_{23}+u_{23}\otimes 1,
$$
\begin{equation}
\Delta u_{13}=1\otimes u_{13}+u_{13}\otimes 1+
u_{12}\otimes u_{23},
\label{a19}
\end{equation}
antipode
\begin{equation}
S(u_{12})=-u_{12}, \quad
S(u_{13})=-u_{13}+u_{12}u_{23}, \quad
S(u_{23})=-u_{23}
\label{a20}
\end{equation}
and standard counit
$\epsilon (u_{12})=\epsilon (u_{13})=\epsilon (u_{23})=0.$

\subsection{\bf Quantum groups
$ SO_v(3;j;\sigma), \; \sigma=(2,1,3) $}

Deformation parameter is transformed by multiplication on
$ J=(\sigma_1,\sigma_3)=(2,3)=j_2. $ 
Commutators, $(v,j)$-orthogonality relations and antipode are
easily obtained from corresponding formulas of   
$ SO_z(3)=SO_v(3;j=1;\sigma_0)$
by interchange of indices 1 and 2 and then by standard reconstruction
of contraction parameters $j.$  In particular, generating matrix is
as follows
\begin{equation}
U(j;\sigma)=\left (\begin{array}{ccc}
u_{22} & j_1u_{21}& j_2u_{23} \\
j_1u_{12} & u_{11} & j_1j_2u_{13} \\
j_2u_{32} & j_1j_2u_{31} & u_{33} \end{array} \right ),
\label{a21}
\end{equation}
Commutation relations of independent generators are
$$
j_1^2[u_{21},u_{13}]=i\frac{1}{j_2}{\rm \sinh } (j_2v) u_{11}(u_{22}-u_{33}),
$$
$$
[u_{23},u_{13}]=u_{13}\left\{\frac{1}{j_2}({\rm \cosh } j_2v-1)u_{23}
-i\frac{1}{j_2}{\rm \sinh } (j_2v) u_{33}\right\},
$$
\begin{equation}
[u_{21},u_{23}]=\left\{\frac{1}{j_2}({\rm \cosh } j_2v-1)u_{23}+
i\frac{1}{j_2}{\rm \sinh } (j_2v) u_{22}\right\}u_{21}.
\label{a22}
\end{equation}
Antipode is easily obtained by the transformations of (\ref{a6})
$$
S(u_{21})=u_{12}{\rm \cosh }(j_2\frac{v}{2})+ij_2^2u_{13}
\frac{1}{j_2}{\rm \sinh }(j_2\frac{v}{2}), \;
$$
$$
S(u_{12})=u_{21}{\rm \cosh }(j_2\frac{v}{2})+ij_2^2u_{31}
\frac{1}{j_2}{\rm \sinh }(j_2\frac{v}{2}),
$$
$$
S(u_{13})=u_{31}{\rm \cosh }(j_2\frac{v}{2})-iu_{21}\frac{1}{j_2}
{\rm \sinh }(j_2\frac{v}{2}),\;
$$
$$
S(u_{31})=u_{13}{\rm \cosh }(j_2\frac{v}{2})-iu_{12}\frac{1}{j_2}
{\rm \sinh }(j_2\frac{v}{2}),
$$
$$
S(u_{23})=u_{32}{\rm \cosh }^2(j_2\frac{v}{2})+u_{23}{\rm \sinh }^2(j_2\frac{v}{2})+
i\frac{1}{2}(u_{33}-u_{22})\frac{1}{j_2}{\rm \sinh }(j_2v),
$$
$$
S(u_{32})=u_{23}{\rm \cosh }^2(j_2\frac{v}{2})+u_{32}{\rm \sinh }^2(j_2\frac{v}{2})+
i\frac{1}{2}(u_{33}-u_{22})\frac{1}{j_2}{\rm \sinh }(j_2v),
$$
$$
S(u_{22})=u_{22}{\rm \cosh }^2(j_2\frac{v}{2}) -u_{33}{\rm \sinh }^2(j_2\frac{v}{2})
+\frac{i}{2}(u_{23}+u_{32})j_2{\rm \sinh }(j_2v),
$$
$$
S(u_{33})=u_{33}{\rm \cosh }^2(j_2\frac{v}{2}) -u_{22}{\rm \sinh }^2(j_2\frac{v}{2})
-\frac{i}{2}(u_{23}+u_{32})j_2{\rm \sinh }(j_2v),
$$
\begin{equation}
\quad S(u_{11})=u_{11}.
\label{a6.1}
\end{equation}
Coproduct and counit are not changed and are given by (\ref{a5}),
which correspond to identical permutation $\sigma_0.$

Contraction  $j_1=\iota_1 $ left deformation parameter fixed
since $J=j_2=1$ and gives {\bf new quantum Euclid group}
$ E_z(2)=SO_z(3;\iota_1,1;\sigma)$  with the matrix
\begin{equation}
U(\iota_1;\sigma)=\left ( \begin{array}{ccc}
\cos\varphi &\iota_1u_{21} & \sin\varphi \\
\iota_1u_{12} &1 & \iota_1u_{13} \\
-\sin\varphi & \iota_1u_{31} & \cos\varphi \end{array} \right )
\sim\left (\begin{array}{ccc}
\cdot & \circ &  \cdot \\
      & \cdot & \circ  \\
 &  & \cdot \end{array} \right ),
\label{a6.1.a}
\end{equation}
where  generators are
$$
u_{11}=1, \;  u_{22}=u_{33}=\cos \varphi, \; u_{23}=-u_{32}=\sin\varphi, 
$$
$$
u_{12}\cos(\varphi-i\frac{v}{2})=-(u_{21} + u_{13}\sin(\varphi-i\frac{v}{2})), 
$$
\begin{equation}
u_{31}\cos(\varphi-i\frac{v}{2})=-(u_{13}+u_{21}\sin(\varphi-i\frac{v}{2})),
\label{a8ad}
\end{equation}
and  the following commutation relations
 $$
  [u_{21},u_{13}]=0,\quad
 [u_{13},\sin\varphi]=2i{\rm \sinh }\frac{z}{2}u_{13}\cos(\varphi-i\frac{z}{2}),
 $$
\begin{equation}
[u_{21},\sin\varphi]=2i{\rm \sinh }\frac{z}{2}\cos(\varphi+i\frac{z}{2})u_{21}
\label{a23}
\end{equation}
are holds. Antipode is given by
$$
S(u_{21})=u_{12}{\rm \cosh }\frac{z}{2}+iu_{13}{\rm \sinh }\frac{z}{2}, \;\;
S(u_{13})=u_{31}{\rm \cosh }\frac{z}{2}-iu_{21}{\rm \sinh }\frac{z}{2},\; \;
$$
\begin{equation}
S(\varphi)=-\varphi,
\label{a6.121e}
\end{equation}
and coproduct is in the form
$$
 \Delta u_{13}= 1\otimes u_{13}+u_{13}\otimes \cos\varphi+
  u_{12}\otimes \sin\varphi,
 $$
\begin{equation}
 \Delta u_{21}=\cos\varphi\otimes u_{21}+u_{21}\otimes 1+
\sin\varphi\otimes u_{31}, \; \;
\Delta \varphi= 1 \otimes \varphi+ \varphi\otimes 1.
\label{a6.121gg}
\end{equation}

{\bf Quantum Newton group } $N_v(2)=SO_v(3;1,\iota_2;\sigma),\; J=\iota_2$
is described by relations
$  u_{33}=1, \;  u_{11}=u_{22}=\cos\psi, \, u_{21}=\sin\psi= -u_{12}, $
i.e. the generating matrix is in the form
\begin{equation}
U(\iota_2;\sigma)=\left ( \begin{array}{ccc}
\cos\psi &\sin\psi & \iota_2u_{23} \\
-\sin\psi &\cos\psi & \iota_2u_{13} \\
\iota_2u_{32} & \iota_2u_{31} & 1 \end{array} \right )
\sim\left (\begin{array}{ccc}
\cdot & \cdot & \bullet \\
 & \cdot & \bullet \\
 &  & \cdot \end{array} \right ),
\label{a6.23g}
\end{equation}
where
$$
u_{31}=-u_{13}\cos\psi - u_{23}\sin\psi -i\frac{v}{2}\sin\psi,\;\;
$$
\begin{equation}
u_{32}=-u_{23}\cos\psi  +u_{13}\sin\psi +i\frac{v}{2}(1-\cos\psi),
\label{a8hg}
\end{equation}
and commutation relations
 $$
 [\sin\psi,u_{13}]=iv\cos\psi(\cos\psi-1),
 $$
\begin{equation}
 [\sin\psi,u_{23}]=iv\cos\psi\sin\psi, \quad
 [u_{23},u_{13}]=-ivu_{13}
\label{a24}
\end{equation}
are holds for independent generators. Antipode is given by
$$ 
S(u_{13})=-u_{13}\cos\psi-u_{23}\sin\psi-iv\sin\psi,\quad
S(\psi)=-\psi,
$$
\begin{equation}
S(u_{23})=-u_{23}\cos\psi+ u_{13}\sin\psi+iv(1-\cos\psi),\; \;
\label{a6.121n}
\end{equation}
and coproduct is
$$
\Delta \psi=1\otimes\psi+\psi\otimes1, \;\;
\Delta u_{23}= u_{23}\otimes 1+\cos\psi\otimes u_{23}+\sin\psi\otimes u_{13},
$$
\begin{equation}
\Delta u_{13}=u_{13}\otimes 1+\cos\psi\otimes u_{13}-\sin\psi\otimes u_{23}.
\label{a6.121nn}
\end{equation}

Generating matrices (\ref{a6.23g}) and (\ref{a12}) are equal 
from the viewpoint of nilpotent units distribution, while formulae
(\ref{a8hg})--(\ref{a6.121nn}) pass to (\ref{a13})--(\ref{a16})
under substitution   $u_{13}$ on $u_{23}$ and $u_{23}$ on $u_{13}.$
Thus, both quantum groups  are isomorphic 
$N_v(2) \simeq N_v^0(2) \simeq  E_v^0(2). $

For {\bf quantum Galilei group }
$G_v(2)=SO_v(3;\iota_1,\iota_2;\sigma),\; J=\iota_2 $
it follows from  $(v,j)$-orthogonality relations that
 $u_{11}=u_{22}=u_{33}=1 $ 
and generating matrix takes the form
\begin{equation}
U(\iota;\sigma)=\left ( \begin{array}{ccc}
1 & \iota_1u_{21} & \iota_2u_{23} \\
-\iota_1u_{12} & 1 & \iota_1\iota_2u_{13} \\
-\iota_2u_{23} & \iota_1\iota_2u_{31} & 1 \end{array} \right )
\sim\left (\begin{array}{ccc}
\cdot & \circ &  \bullet \\
 & \cdot & \times \\
 &  & \cdot \end{array} \right ),
\label{a8.21g}
\end{equation}
where  $u_{31}=-u_{13}-u_{21}u_{23}+i\frac{v}{2}u_{21}, $ 
commutation relations are
\begin{equation}
 [u_{21},u_{13}]=0, \quad
 [u_{23},u_{13}]=-ivu_{13}, \quad
 [u_{21},u_{23}]=ivu_{21},
\label{a25}
\end{equation}
antipode may be written as
\begin{equation}
S(u_{21})=-u_{21},\quad
S(u_{23})=-u_{23},\quad
S(u_{13})=-u_{13}-u_{21}u_{23}, 
\label{a6.121g}
\end{equation}
and coproduct is
$$
\Delta u_{21}=1\otimes u_{21}+u_{21}\otimes 1, \;\;
\Delta u_{23}=1\otimes u_{23}+u_{23}\otimes 1,
$$
\begin{equation}
\Delta u_{13}=1\otimes u_{13}+u_{13}\otimes 1+u_{21}\otimes u_{23}.
\label{a6.121gn}
\end{equation}

Let us stress that  $G_v(2)$ is not isomorphic to $G_v^0(2),$
in spite of the fact that both matrices  (\ref{a8.21g}), (\ref{a17})
are equivalent from the viewpoint of nilpotent units distribution,
but deformation parameters are transformed in a different ways,
namely, with multipliers $ J=\iota_2 $ and  $ J=\iota_1\iota_2 $
respectively. Therefore commutation relations (\ref{a18}), (\ref{a25}),
antipodes (\ref{a20}), (\ref{a6.121g})  and counits are passed in each other
under substitution  $u_{13}$ on $ u_{23}$ and vice versa, but in coproduct
(\ref{a19}) $\Delta (u_{13})$ is not passed in  $ \Delta(u_{23})$ from
(\ref{a6.121gn}).

\subsection{\bf Quantum groups $ SO_v(3;j;\sigma), \; \sigma=(1,3,2) $}

Deformation parameter is multiplied by
$ J=(\sigma_1,\sigma_3)=(1,2)=j_1. $ 
Commutators, $(v,j)$-orthogonality relations and antipode are
easily obtained from corresponding formulas of   
$ SO_z(3)=SO_v(3;1,1;\sigma_0)$
by interchange of indices 2 and 3 and then by standard reconstruction
of contraction parameters $j.$  In particular, generating matrix is
as follows
\begin{equation}
U(j;\sigma)=\left (\begin{array}{ccc}
u_{11} & j_1j_2u_{13}& j_1u_{12} \\
j_1j_2u_{31} & u_{33} & j_2u_{32} \\
j_1u_{21} & j_2u_{23} & u_{22} \end{array} \right ).
\label{a26}
\end{equation}
For $j_1=\iota_1 $ {\bf quantum Euclid group}
$ \tilde{E}_v(2)=SO_v(3;\iota_1,1;\sigma)$ is obtained
with generators
\begin{equation}
U(\iota_1;\sigma)=\left ( \begin{array}{ccc}
1 & \iota_1u_{13} & \iota_1u_{12} \\
\iota_1u_{31} & \cos\varphi & \sin\varphi \\
\iota_1u_{21} & -\sin\varphi & \cos\varphi \end{array} \right )
\sim\left (\begin{array}{ccc}
\cdot & \circ & \circ \\
      & \cdot &  \cdot \\
 &  & \cdot \end{array} \right ).
\label{a6.253}
\end{equation}
As far as the generating matrix (\ref{a6.253}) is equal to (\ref{a7}),
then  $ \tilde{E}_v(2)$ is isomorphic with  $E^0_v(2)$ and therefore
do not present a new quantum group.

{\bf Quantum Newton group}
$ \tilde{N}_z(2)=SO_z(3;1,\iota_2;\sigma)$ is described by
untouched deformation parameter  $z,$ generators
$  u_{33}=1,\;  u_{11}=u_{22}=\cos\psi, \; u_{12}=\sin\psi=-u_{21}, $
which are arranged in matrix form
\begin{equation}
U(\iota_2;\sigma)=\left ( \begin{array}{ccc}
\cos\psi & \iota_2u_{13} & \sin\psi \\
\iota_2u_{31} & 1 & \iota_2u_{32} \\
-\sin\psi & \iota_2u_{23} & \cos\psi \end{array} \right )
\sim\left (\begin{array}{ccc}
\cdot & \bullet  & \cdot \\
 & \cdot & \bullet \\
 &  & \cdot \end{array} \right ).
\label{a8-dd}
\end{equation}
This quantum group  as Hopf algebra is isomorphic to quantum Euclid group
$E_z(2)$ with untouched deformation parameter $(J=1),$
since the generating matrix (\ref{a8-dd}) is equal to (\ref{a6.1.a}),
if instead of $\iota_{2}$ put $\iota_{1}.$
Finally,  {\bf quantum Galilei group}
$ \tilde{G}_v(2)=SO_v(3;\iota_1,\iota_2;\sigma)$
is characterized by   $J=\iota_1,$ 
diagonal generators are equal to one
 $u_{11}=u_{22}=u_{33}=1, $  and generating matrix is as follows
\begin{equation}
U(\iota;\sigma)=\left ( \begin{array}{ccc}
1 & \iota_1\iota_2u_{13} & \iota_1u_{12} \\
\iota_1\iota_2u_{31} & 1 & \iota_2u_{32} \\
-\iota_1u_{12} & -\iota_2u_{32} & 1 \end{array} \right )
\sim \left (\begin{array}{ccc}
\cdot & \times  &    \circ \\
 & \cdot & \bullet \\
 &  & \cdot \end{array} \right ).
\label{a9.32}
\end{equation}
 The nilpotent parameters distribution of (\ref{a9.32}) pass
in (\ref{a8.21g}) under   exchange  $ \iota_1 $  and $ \iota_2, $  
and  simultaneous reflection with respect to secondary diagonal.
Therefore,  $\tilde{G}_v(2)$ is isomorphic to $G_v(2).$
Thus the permutation $ \sigma=(1,3,2) $ do not lead to new
contracted quantum groups.

\section{Quantum groups $ SO_v(4;j;\sigma)$} 

In this section all nonisomorphic contractions of $SO_q(4)$ are enumerated.
Deformation parameter is multiplied on
$ J=(\sigma_1,\sigma_4) \cup (\sigma_2,\sigma_3),$ 
which is equal to 
$ J=j_1j_2j_3 $ for permutation  $ \sigma_0 =(1,2,3,4)$ and
$ J=j_1j_3 $ for  $ \sigma' =(1,3,4,2).$ 
There are not other values of $J.$ Above-mentioned values of $J$
correspond to nonisomorphic  on the equal parameter number contracted
quantum groups which have nonequivalent generating matrices for
permutations $ \sigma_0 $ and $ \sigma'.$

{\bf One-dimensional contractions.} 
For $j_1=\iota_1, \;  J=\iota_1 $ quantum Euclid group 
$E_v(3)=SO_v(4;\iota_1;\sigma_0)$ is obtained.
For $j_2=\iota_2$ there are two nonisomorphic quantum Newton groups:
$N_v(3)=SO_v(4;\iota_2;\sigma_0), \; J=\iota_2 $ and
$N_z(3)=SO_z(4;\iota_2;\sigma')$  with  $ J=1. $
{\bf Two-dimensional contractions.} 
For $j_1=\iota_1, j_2=\iota_2  $
two nonisomorphic quantum Galilei groups:
$G_v(3)=SO_v(4;\iota_1,\iota_2;\sigma_0), \;  J=\iota_1\iota_2 $  and
$G_w(3)=SO_w(4;\iota_1,\iota_2;\sigma'), \;  J=\iota_1 $ are obtained.
Contractions  $j_1=\iota_1, j_3=\iota_3  $ gives in result
 quantum  groups $ SO_v(4;\iota_1,\iota_3;\sigma_0),$ $ J=\iota_1\iota_3, $
which has not special name.
Under maximal {\bf three-dimensional contractions} 
$j_1=\iota_1, j_2=\iota_2,  j_3=\iota_3  $
two nonisomorphic quantum flag groups:
$F_v(4)= SO_v(4;\iota;\sigma_0), \;  J=\iota_1\iota_2\iota_3 $ and
$F_w(4)= SO_w(4;\iota;\sigma'), \;  J=\iota_1\iota_3 $ are obtained.

$$
 E_v(3) \sim
 \left (\begin{array}{cccc}
\cdot & \circ & \circ & \circ \\
      & \cdot & \cdot & \cdot \\
      &       & \cdot & \cdot \\
      &       &       & \cdot 
\end{array} \right), \quad
N_v(3)\sim\left (\begin{array}{cccc}
\cdot & \cdot & \bullet & \bullet \\
      & \cdot & \bullet & \bullet \\
      &       & \cdot & \cdot \\
      &       &       & \cdot 
\end{array} \right),
$$
$$
N_z(3)\sim\left (\begin{array}{cccc}
\cdot & \bullet & \bullet & \cdot \\
      & \cdot   & \cdot  & \bullet \\
      &         & \cdot  & \bullet \\
      &         &        & \cdot 
\end{array} \right), \quad
G_v(3) \sim\left (\begin{array}{cccc}
\cdot & \circ   & \times  & \times  \\
      & \cdot   & \bullet & \bullet \\
      &         & \cdot   & \cdot   \\
      &         &         & \cdot   
\end{array} \right),
$$
$$
G_w(3) \sim\left (\begin{array}{cccc}
\cdot & \times  & \times  & \circ    \\
      & \cdot   & \cdot   & \bullet  \\
      &         & \cdot   & \bullet  \\
      &         &         & \cdot   
\end{array} \right), \quad
F_v(4)\sim\left (\begin{array}{cccc}
\cdot & \circ  &  \times  & \otimes     \\
      & \cdot   & \bullet & \diamondsuit  \\
      &         & \cdot   & \star        \\
      &         &         & \cdot       
\end{array} \right), 
$$
$$
SO_v(4;\iota_1,\iota_3;\sigma_0)\sim\left (\begin{array}{cccc}
\cdot & \circ  &  \circ  & \triangle    \\
      & \cdot   & \cdot   & \star        \\
      &         & \cdot   & \star        \\
      &         &         & \cdot       
\end{array} \right), \quad
F_w(4)\sim\left (\begin{array}{cccc}
\cdot & \times  & \otimes &  \circ   \\
      & \cdot   &  \star   &  \bullet     \\
      &         &  \cdot  & \diamondsuit  \\
      &         &         & \cdot         
\end{array} \right),
$$
where
$\triangle = \iota_3,\;\star =\iota_1\iota_3, \;
\diamondsuit=\iota_2\iota_3,\; \otimes = \iota_1\iota_2\iota_3. $

Thus, for  quantum case there are eight different contracted groups
while for classical group $SO(4)$ there are only five nonisomorphic
contracted Caley--Klein groups.

\section{Quantum groups $ SO_v(5;j;\sigma)$}

Deformation parameter is multiplied on
$ J=(\sigma_1,\sigma_5) \cup (\sigma_2,\sigma_4),$ 
which is equal to 
$ J=j_1j_2j_3j_4 $ for permutation $ \sigma_0 =(1,2,3,4,5),$  
equal to $ J=j_1j_2j_3 $ for permutation $ \sigma^1 =(1,2,5,3,4),$ 
equal to $ J=j_1j_2j_4 $ for permutation $ \sigma^2 =(1,4,2,5,3),$ 
equal to $ J=j_1j_3 $ for permutation $ \sigma^3 =(1,3,5,4,2),$ 
equal to $ J=j_1j_4 $ for permutation $ \sigma^4 =(1,4,3,5,2),$ 
equal to $ J=j_2j_4 $ for permutation $ \sigma^5 =(2,4,1,5,3),$ 
equal to $ J=j_1j_3j_4 $ for permutation $ \sigma^6 =(1,3,4,5,2),$ 
equal to $ J=j_2j_3j_4 $ for permutation $ \sigma^7 =(2,3,1,4,5).$ 

If  contractions only  on parameters  $ j_1, j_2 $ are considered,
then two quantum Euclid groups:
$E_v(4)=SO_v(4;\iota_1;\sigma_0), \; J=\iota_1 $ and
$E_z(4)=SO_z(4;\iota_1;\sigma^5), \; J=1 $
with distribution of nilpotent parameters in the form
$$
E_v(4)\sim\left (\begin{array}{ccccc}
\cdot & \circ  &  \circ  & \circ & \circ    \\
      & \cdot  & \cdot   & \cdot & \cdot    \\
      &        & \cdot   & \cdot & \cdot    \\
      &        &         & \cdot & \cdot    \\  
      &        &         &       & \cdot    
\end{array} \right), \qquad
E_z(4)\sim\left (\begin{array}{ccccc}
\cdot & \circ  &  \cdot  & \cdot & \cdot    \\
      & \cdot  & \circ   & \circ & \circ    \\
      &        & \cdot   & \cdot & \cdot    \\
      &        &         & \cdot & \cdot    \\  
      &        &         &       & \cdot    
\end{array} \right), 
$$
two quantum Newton groups:
$N_v(4)=SO_v(4;\iota_2;\sigma_0),$ $ J=\iota_2 $ and
$N_z(4)=SO_z(4;\iota_2;\sigma^3),$ $ J=1 $
with generating matrices
$$
N_v(4)\sim\left (\begin{array}{ccccc}
\cdot & \cdot  &  \bullet  & \bullet & \bullet    \\
      & \cdot  & \bullet   & \bullet & \bullet    \\
      &        & \cdot   & \cdot & \cdot    \\
      &        &         & \cdot & \cdot    \\  
      &        &         &       & \cdot    
\end{array} \right), \qquad
N_z(4)\sim\left (\begin{array}{ccccc}
\cdot & \bullet  &  \cdot  & \bullet & \bullet  \\
      & \cdot  & \bullet   & \cdot & \cdot    \\
      &        & \cdot   & \bullet & \bullet    \\
      &        &         & \cdot & \cdot    \\  
      &        &         &       & \cdot    
\end{array} \right), 
$$
and two quantum Galilei groups:
$G_v(4)=SO_v(4;\iota_1\iota_2;\sigma_0),$ $ J=\iota_1\iota_2 $ and
$G_z(4)=SO_v(4;\iota_1\iota_2;\sigma^3),$ $ J=\iota_1 $
with generating matrices
$$
G_v(4)\sim\left (\begin{array}{ccccc}
\cdot & \circ  &  \times  & \times & \times    \\
      & \cdot  & \bullet   & \bullet & \bullet    \\
      &        & \cdot   & \cdot & \cdot    \\
      &        &         & \cdot & \cdot    \\  
      &        &         &       & \cdot    
\end{array} \right), \qquad
G_z(4)\sim\left (\begin{array}{ccccc}
\cdot & \times &  \circ  & \times & \times    \\
      & \cdot  & \bullet   & \cdot & \cdot    \\
      &        & \cdot   & \bullet & \bullet    \\
      &        &         & \cdot & \cdot    \\  
      &        &         &       & \cdot    
\end{array} \right). 
$$
As compared with the case  $N=3$ two quantum Newton groups are added.

In all discussed examples for  $N=3,4,5$ the number of nonisomorphic
quantum analogues of the corresponding classical groups is equal two.
It may be  think that this number for any contractions do not exceed two.
But this is not so. The number of nonisomorphic quantum analogues 
of the  classical Caley--Klein groups is increased when the number of
nilpotent valued contraction parameters is increased. For example,
under maximal contraction $j_k=\iota_k, k=1,\ldots,4  $  five
quantum analogues of the flag group $ F(5)=SO(5;\iota)$ are obtained, 
namely:
$F_v(5)= SO_v(5;\iota;\sigma_0), \;  J=\iota_1\iota_2\iota_3\iota_4; $
$F_{v_1}(5)= SO_{v_1}(5;\iota;\sigma^1), \;  J=\iota_1\iota_2\iota_3; $
$F_{v_2}(5)= SO_{v_2}(5;\iota;\sigma^2), \;  J=\iota_1\iota_2\iota_4; $
$F_{v_3}(5)= SO_{v_3}(5;\iota;\sigma^3), \;  J=\iota_1\iota_3; $
$F_{v_4}(5)= SO_{v_4}(5;\iota;\sigma^4), \;  J=\iota_1\iota_4. $
All they have generating matrices with nonequivalent distributions
of nilpotent parameters.

\section*{Acknowledgments}

N.G. is grateful to P.P.Kulish and V.O.Tarasov for fruitful  discussions.

\vspace{5mm}

\hfill Appendix 1
\begin{center}
    {\large \bf $R$-matrix of quantum group $ SO_q(N)$
             in Cartesian basis }
\end{center}

$$
\tilde{R}_q=(D \otimes D) R  ( D \otimes D)^{-1}  =
$$
$$
=I+ \frac{1}{2}(q-1)(1-q^{-1})\sum_{k=1 \atop k\neq k'}^N
        ( e_{kk} \otimes  e_{kk}
        + e_{kk} \otimes  e_{k'k'})
+ \frac{\lambda}{2} \sum_{k=1 \atop k\neq k'}^N
        (   e_{k'k} \otimes  e_{kk'}
         -  e_{k'k} \otimes  e_{k'k} )
$$
$$
+ \frac{\lambda}{2} \sum_{k=1}^n 
(        e_{k',n+1} \otimes  e_{n+1,k'} 
      -i e_{k',n+1} \otimes  e_{n+1,k} 
      +i e_{k,n+1} \otimes  e_{n+1,k'} 
      +  e_{k,n+1} \otimes  e_{n+1,k} 
$$
$$
   +  e_{n+1,k} \otimes  e_{k,n+1} 
   +i e_{n+1,k} \otimes  e_{k',n+1} 
   -i e_{n+1,k'} \otimes  e_{k,n+1} 
   +  e_{n+1,k'} \otimes  e_{k',n+1} ) 
$$
$$
- \frac{\lambda}{2} \sum_{k=1 }^n q^{-\rho_k} 
( -i e_{k',n+1} \otimes   e_{k,n+1} 
  +  e_{k',n+1} \otimes   e_{k',n+1} 
  +  e_{k,n+1} \otimes   e_{k,n+1} 
  +i e_{k,n+1} \otimes   e_{k',n+1}    
$$
$$
   +i  e_{n+1,k} \otimes   e_{n+1,k'}  
   +   e_{n+1,k} \otimes   e_{n+1,k} 
   +   e_{n+1,k'} \otimes  e_{n+1,k'} 
   -i  e_{n+1,k'} \otimes  e_{n+1,k}) 
$$
$$
+ \frac{\lambda}{4} \sum_{k,p=1 \atop k > p,\; k,p\neq n+1}^N 
(     e_{kp} \otimes  e_{pk}  
   +  e_{kp} \otimes  e_{p'k'}  
   +i e_{kp} \otimes  e_{p'k} 
   -i e_{kp} \otimes  e_{pk'} 
$$
$$
   +  e_{k'p'} \otimes  e_{pk}  
   +  e_{k'p'} \otimes  e_{p'k'}  
   +i e_{k'p'} \otimes  e_{p'k} 
   -i e_{k'p'} \otimes  e_{pk'} 
$$
$$
   +i e_{k'p} \otimes  e_{pk}  
   +i e_{k'p} \otimes  e_{p'k'}  
   -  e_{k'p} \otimes  e_{p'k} 
   +  e_{k'p} \otimes  e_{pk'} 
$$
$$
   -i e_{kp'} \otimes  e_{pk}  
   -i e_{kp'} \otimes  e_{p'k'}  
   +  e_{kp'} \otimes  e_{p'k} 
   -  e_{kp'} \otimes  e_{pk'}) 
$$
$$
-\frac{\lambda}{4} \sum_{k,p=1 \atop k>p,\; k,p \neq n+1}^N q^{\rho_k-\rho_p} 
(     e_{kp} \otimes   e_{k'p'}  
   +  e_{kp} \otimes   e_{kp}  
   +i e_{kp} \otimes   e_{kp'} 
   -i e_{kp} \otimes   e_{k'p} 
$$
$$
   +  e_{k'p'} \otimes   e_{k'p'}  
   +  e_{k'p'} \otimes   e_{kp}  
   +i e_{k'p'} \otimes   e_{kp'} 
   -i e_{k'p'} \otimes   e_{k'p} 
$$
$$
   +i e_{k'p} \otimes   e_{k'p'}  
   +i e_{k'p} \otimes   e_{kp}  
   -  e_{k'p} \otimes   e_{kp'} 
   +  e_{k'p} \otimes   e_{k'p} 
$$
$$
   -i e_{kp'} \otimes  e_{k'p'} 
   -i e_{kp'} \otimes  e_{kp} 
   +  e_{kp'} \otimes  e_{kp'} 
   -  e_{kp'} \otimes  e_{k'p}), \quad  \lambda =q-q^{-1}.
$$

\vspace{15mm}

\hfill Appendix 2
\begin{center}
    {\large \bf Antipode of quantum group  $ SO_v(N;j;\sigma)$
              in Cartesian basis }
\end{center}
Antipode of Cartesian generators of quantum group
$ SO_v(N;j;\sigma),$ $ N=2n+1 $ is obtained by formula
$$
S(U(j;\sigma))={\tilde C}_v(j)U^t(j;\sigma){\tilde C}^{-1}_v(j)
$$
with the help of matrix    ${\tilde C}_v(j)=D^{-1}C_v(j)(D^t)^{-1} $
in the following form
$$
S(u_{\sigma_k\sigma_{n+1}})=u_{\sigma_{n+1}\sigma_{k}}\cosh (Jv{\rho}_k)
+iu_{\sigma_{n+1}\sigma_{k'}}
\frac{(\sigma_{k'},\sigma_{n+1})}{(\sigma_k,\sigma_{n+1})}
\sinh(Jv\rho_k),
$$
$$
S(u_{\sigma_{n+1}\sigma_k})=u_{\sigma_k\sigma_{n+1}}\cosh(Jv{\rho}_k)
+iu_{\sigma_{k'}\sigma_{n+1}}
\frac{(\sigma_{k'},\sigma_{n+1})}{(\sigma_k,\sigma_{n+1})}
\sinh(Jv{\rho}_k),
$$
$$
S(u_{\sigma_{n+1+k}\sigma_{n+1}})=u_{\sigma_{n+1}\sigma_{n+1+k}}
\cosh(Jv{\rho}_{n+1-k})-
$$
$$
-iu_{\sigma_{n+1}\sigma_{n+1-k}}
\frac{(\sigma_{n+1-k},\sigma_{n+1})}{(\sigma_{n+1+k},\sigma_{n+1})}
\sinh(Jv{\rho}_{n+1-k}),
$$
$$
S(u_{\sigma_{n+1}\sigma_{n+1+k}})=u_{\sigma_{n+1+k}\sigma_{n+1}}
\cosh(Jv{\rho}_{n+1-k})-
$$
$$
-iu_{\sigma_{n+1-k}\sigma_{n+1}}
\frac{(\sigma_{n+1-k},\sigma_{n+1})}{(\sigma_{n+1+k},\sigma_{n+1})}
\sinh(Jv\rho_{n+1-k}),
$$
$$
S(u_{\sigma_k\sigma_p})=u_{\sigma_p\sigma_k}
\cosh(Jv{\rho}_{k})\cosh(Jv{\rho}_p)-
$$
$$
-u_{\sigma_{p'}\sigma_{k'}}
\frac{(\sigma_{k'},\sigma_{p'})}{(\sigma_k,\sigma_p)}
\sinh(Jv{\rho}_{k})\sinh(Jv{\rho}_{p})+
$$
$$
+i\left ( u_{\sigma_p\sigma_{k'}}
\frac{(\sigma_{k'},\sigma_p)}{(\sigma_k,\sigma_p)}
\sinh(Jv{\rho}_{k})\cosh(Jv{\rho}_{p})+ \right.
$$
$$
+ \left. u_{\sigma_{p'}\sigma_k}
\frac{(\sigma_k,\sigma_{p'})}{(\sigma_k,\sigma_p)}
\cosh(Jv{\rho}_{k})\sinh(Jv{\rho}_{p}) \right ),
$$
$$ 
S(u_{\sigma_k\sigma_{n+1+p}})=u_{\sigma_{n+1+p}\sigma_k}
\cosh(Jv{\rho}_{k})\cosh(Jv{\rho}_{{n+1-p}})+
$$
$$
+u_{\sigma_{n+1-p}\sigma_{k'}}
\frac{(\sigma_{k'},\sigma_{n+1-p})}{(\sigma_k,\sigma_{n+1+p})}
\sinh(Jv{\rho}_{k})\sinh(Jv{\rho}_{{n+1-p}})+
$$
$$
+i\left ( u_{\sigma_{n+1+p}\sigma_{k'}}
\frac{(\sigma_{k'},\sigma_{n+1+p})}{(\sigma_k,\sigma_{n+1+p})}
\sinh(Jv{\rho}_{k})\cosh(Jv{\rho}_{{n+1-p}})- \right.
$$
$$
\left. -u_{\sigma_{n+1-p}\sigma_k}
\frac{(\sigma_k,\sigma_{n+1-p})}{(\sigma_k,\sigma_{n+1+p})}
\cosh(Jv{\rho}_{k})\sinh(Jv{\rho}_{{n+1-p}}) \right ),
$$
$$
S(u_{\sigma_{n+1+k}\sigma_p})=u_{\sigma_p\sigma_{n+1+k}}
\cosh(Jv{\rho}_{{n+1-k}})\cosh(Jv{\rho}_{p})+
$$
$$
+u_{\sigma_{p'}\sigma_{n+1-k}}
\frac{(\sigma_{n+1-k},\sigma_{p'})}{(\sigma_{n+1+k},\sigma_p)}
\sinh(Jv{\rho}_{{n+1-k}})\sinh(Jv{\rho}_{p})+
$$
$$
+i\left ( u_{\sigma_{p'}\sigma_{n+1+k}}
\frac{(\sigma_{n+1+k},\sigma_{p'})}{(\sigma_{n+1+k},\sigma_p)}
\cosh(Jv{\rho}_{{n+1-k}})\sinh(Jv{\rho}_{p})- \right.
$$
$$
\left. -u_{\sigma_p\sigma_{n+1-k}}
\frac{(\sigma_{n+1-k},\sigma_p)}{(\sigma_{n+1+k},\sigma_p)}
\sinh(Jv{\rho}_{{n+1-k}})\cosh(Jv{\rho}_{p}) \right ),
$$
$$
S(u_{\sigma_{n+1+k}\sigma_{n+1+p}})=
u_{\sigma_{n+1+p}\sigma_{n+1+k}}
\cosh(Jv{\rho}_{{n+1-k}})\cosh(Jv{\rho}_{{n+1-p}})-
$$
$$
-u_{\sigma_{n+1-p}\sigma_{n+1-k}}
\frac{(\sigma_{n+1-k},\sigma_{n+1-p})}{(\sigma_{n+1+k},\sigma_{n+1+p})}
\sinh(Jv{\rho}_{{n+1-k}})\sinh(Jv{\rho}_{{n+1-p}})-
$$
$$
-i\left ( u_{\sigma_{n+1-p}\sigma_{n+1+k}}
\frac{(\sigma_{n+1+k},\sigma_{n+1-p})}{(\sigma_{n+1+k},\sigma_{n+1+p})}
\cosh(Jv{\rho}_{{n+1-k}})\sinh(Jv{\rho}_{{n+1-p}})+ \right.
$$
\begin{equation}
\left. +u_{\sigma_{n+1+p}\sigma_{n+1-k}}
\frac{(\sigma_{n+1-k},\sigma_{n+1+p})}{(\sigma_{n+1+k},\sigma_{n+1+p})}
\sinh(Jv{\rho}_{{n+1-k}})\cosh(Jv{\rho}_{{n+1-p}}) \right ),  
   \label{a-6}
\end{equation}
where  $ k,p=1, \ldots,n. $ 
Antipode of  quantum group $ SO_v(N;j;\sigma), \ N=2n $
is given by above-mentioned formulae with the replacement
$n+1$ on $n.$

\vspace{5mm}

\hfill Appendix 3
\begin{center}
    {\large \bf $(v,j)$-orthogonality relations of quantum group
             $ SO_v(N;\sigma;j) $ in Cartesian basis }
\end{center}

Additional relations
$ U(j;\sigma)\tilde{C}_v(j)U^t(j;\sigma)=\tilde{C}_v(j), $ 
where
$\tilde{C}_v(j)=D^{-1}C_v(j)(D^t)^{-1}$
are in the form 
$$ 
u_{\sigma_k\sigma_{n+1}}
u_{\sigma_p\sigma_{n+1}}(\sigma_k,\sigma_{n+1})
(\sigma_p,\sigma_{n+1})+
$$
$$
+\sum^n_{s=1}
\left\{ u_{\sigma_k\sigma_s}u_{\sigma_p\sigma_s}
(\sigma_k,\sigma_s)(\sigma_p,\sigma_s)
\cosh(Jv{\rho}_{{s}})+ \right.
$$
$$
+u_{\sigma_k\sigma_{n+1+s}}u_{\sigma_p\sigma_{n+1+s}}
(\sigma_k,\sigma_{n+1+s})(\sigma_p,\sigma_{n+1+s})
\cosh(Jv{\rho}_{{n+1-s}})+
$$
$$
 +i\left [ u_{\sigma_k\sigma_{n+1-s}}
u_{\sigma_p\sigma_{n+1+s}}(\sigma_k,\sigma_{n+1-s})
(\sigma_p,\sigma_{n+1+s})
\sinh(Jv{\rho}_{{n+1-s}})- \right.
$$
 \begin{equation}
\left. \left.-u_{\sigma_k\sigma_{s'}}
u_{\sigma_p\sigma_s}
(\sigma_k,\sigma_{s'})(\sigma_p,\sigma_s)
\sinh(Jv{\rho}_{{s}}) \right ] \right\}=
\delta_{kp} \cosh(Jv{\rho}_{{k}}),
\label{o1-1}
\end{equation}
$$
 u_{\sigma_{n+1+k}\sigma_{n+1}}
u_{\sigma_{n+1+p}\sigma_{n+1}}
(\sigma_{n+1+k},\sigma_{n+1})(\sigma_{n+1+p},\sigma_{n+1})+
$$
$$
+\sum^n_{s=1}
\left\{ u_{\sigma_{n+1+k}\sigma_s}u_{\sigma_{n+1+p}\sigma_s}
(\sigma_{n+1+k},\sigma_s)(\sigma_{n+1+p},\sigma_s)
\cosh(Jv{\rho}_{{s}})+       \right.
$$
$$
+u_{\sigma_{n+1+k}\sigma_{n+1+s}}u_{\sigma_{n+1+p}\sigma_{n+1+s}}
(\sigma_{n+1+k},\sigma_{n+1+s})(\sigma_{n+1+p},\sigma_{n+1+s})
 \cosh(Jv{\rho}_{{n+1-s}})+
$$
$$
+ i\left [ u_{\sigma_{n+1+k}\sigma_{n+1-s}}
 u_{\sigma_{n+1+p}\sigma_{n+1+s}}
(\sigma_{n+1+k},\sigma_{n+1-s})(\sigma_{n+1+p},\sigma_{n+1+s})
\sinh(Jv{\rho}_{{n+1-s}})- \right.
$$
   \begin{equation}
\left. \left. -u_{\sigma_{n+1+k}\sigma_{s'}}
u_{\sigma_{n+1+p}\sigma_s}
(\sigma_{n+1+k},\sigma_{s'})(\sigma_{n+1+p},\sigma_s)
\sinh(Jv{\rho}_{{s}}) \right ] \right\}=
$$
$$
=\delta_{kp} \cosh(Jv{\rho}_{{n+1-k}}),
\label{o1-2}
\end{equation}
$$
u_{\sigma_k\sigma_{n+1}}u_{\sigma_{n+1+p}\sigma_{n+1}}
(\sigma_k,\sigma_{n+1})(\sigma_{n+1+p},\sigma_{n+1})+
$$
$$
+\sum^n_{s=1}
\left\{ u_{\sigma_k\sigma_s}u_{\sigma_{n+1+p}\sigma_s}
(\sigma_k,\sigma_s)(\sigma_{n+1+p},\sigma_s)
\cosh(Jv{\rho}_{{s}})+ \right.
$$
$$
+u_{\sigma_k\sigma_{n+1+s}}u_{\sigma_{n+1+p}\sigma_{n+1+s}}
(\sigma_k,\sigma_{n+1+s})(\sigma_{n+1+p},\sigma_{n+1+s})
\cosh(Jv{\rho}_{{n+1-s}})+
$$
$$
 +i\left [ u_{\sigma_k\sigma_{n+1-s}}u_{\sigma_{n+1+p}\sigma_{n+1+s}}
(\sigma_k,\sigma_{n+1-s})(\sigma_{n+1+p},\sigma_{n+1+s})
\sinh(Jv{\rho}_{{n+1-s}})- \right.
$$
  \begin{equation}
\left. \left. -u_{\sigma_k\sigma_{s'}}u_{\sigma_{n+1+p}\sigma_s}
(\sigma_k,\sigma_{s'})(\sigma_{n+1+p},\sigma_s)
\sinh(Jv{\rho}_{{s}}) \right ] \right\}=
i\delta_{n+1-k,p} \sinh(Jv{\rho}_{{k}}),
\label{o1-3}
\end{equation}
$$
u_{\sigma_{n+1+k}\sigma_{n+1}}u_{\sigma_p\sigma_{n+1}}
(\sigma_{n+1+k},\sigma_{n+1})(\sigma_p,\sigma_{n+1})+
$$
$$
+\sum^n_{s=1}
\left\{ u_{\sigma_{n+1+k}\sigma_s}u_{\sigma_p\sigma_s}
(\sigma_{n+1+k},\sigma_s)(\sigma_p,\sigma_s)
\cosh(Jv{\rho}_{{s}})+ \right.
$$
$$
+u_{\sigma_{n+1+k}\sigma_{n+1+s}}
u_{\sigma_p\sigma_{n+1+s}}
(\sigma_{n+1+k},\sigma_{n+1+s})(\sigma_p,\sigma_{n+1+s})
\cosh(Jv{\rho}_{{n+1-s}})+
$$
$$
 +i\left [ u_{\sigma_{n+1+k}\sigma_{n+1-s}}
u_{\sigma_p\sigma_{n+1+s}}
(\sigma_{n+1+k},\sigma_{n+1-s})(\sigma_p,\sigma_{n+1+s})
\sinh(Jv{\rho}_{{n+1-s}})-  \right.
$$
  \begin{equation}
\left. \left.-u_{\sigma_{n+1+k}\sigma_{s'}}u_{\sigma_p\sigma_s}
(\sigma_{n+1+k},\sigma_{s'})(\sigma_p,\sigma_s)
\sinh(Jv{\rho}_{{s}}) \right ] \right\}=
-i\delta_{n+1-k,p} \sinh(Jv{\rho}_{{p}}),
\label{o1-4}
\end{equation}
$$
u^2_{\sigma_{n+1}\sigma_{n+1}}
+\sum^n_{k=1}
\left\{ u^2_{\sigma_{n+1}\sigma_k}
(\sigma_{n+1},\sigma_k)^2
\cosh(Jv{\rho}_{{k}})+ \right.
$$
$$
+u^2_{\sigma_{n+1}\sigma_{n+1+k}}
(\sigma_{n+1},\sigma_{n+1+k})^2
\cosh(Jv{\rho}_{{n+1-k}})+
$$
$$
+ i\left [ u_{\sigma_{n+1}\sigma_k}u_{\sigma_{n+1}\sigma_{k'}}
(\sigma_{n+1},\sigma_k)(\sigma_{n+1},\sigma_{k'})
\sinh(Jv{\rho}_k)- \right.
$$
  \begin{equation}
\left. \left.-u_{\sigma_{n+1}\sigma_{n+1+k}}
u_{\sigma_{n+1}\sigma_{n+1-k}}
(\sigma_{n+1},\sigma_{n+1+k})(\sigma_{n+1},\sigma_{n+1-k})
\sinh(Jv{\rho}_{{n+1-k}}) \right ] \right\}=1,
\label{o1-5}
\end{equation}
$$
u_{\sigma_k\sigma_{n+1}}u_{\sigma_{n+1}\sigma_{n+1}}
(\sigma_k,\sigma_{n+1})
+\sum^n_{s=1}
\left\{ u_{\sigma_k\sigma_s}u_{\sigma_{n+1}\sigma_s}
(\sigma_k,\sigma_s)(\sigma_{n+1},\sigma_s)
\cosh(Jv{\rho}_{s})+ \right.
$$
$$
+u_{\sigma_k\sigma_{n+1+s}}
u_{\sigma_{n+1}\sigma_{n+1+s}}
(\sigma_k,\sigma_{n+1+s})(\sigma_{n+1},\sigma_{n+1+s})
\cosh(Jv{\rho}_{{n+1-s}})+
$$
$$
+i\left [ u_{\sigma_k\sigma_s}
u_{\sigma_{n+1}\sigma_{s'}}
(\sigma_k,\sigma_s)(\sigma_{n+1},\sigma_{s'})
\sinh(Jv{\rho}_{{s}})- \right.
$$
   \begin{equation}
\left. \left. -u_{\sigma_k\sigma_{n+1+s}}
u_{\sigma_{n+1}\sigma_{n+1-s}}
(\sigma_k,\sigma_{n+1+s})(\sigma_{n+1},\sigma_{n+1-s})
\sinh(Jv{\rho}_{{n+1-s}}) \right ] \right\}=0,
\label{o1-6}
\end{equation}
$$
u_{\sigma_{n+1}\sigma_{n+1}}
u_{\sigma_k\sigma_{n+1}}
(\sigma_k,\sigma_{n+1})
+\sum^n_{s=1}
\left\{ u_{\sigma_{n+1}\sigma_s}u_{\sigma_k\sigma_s}
(\sigma_{n+1},\sigma_s)(\sigma_k,\sigma_s)
\cosh(Jv{\rho}_{{s}})+ \right.
$$
$$
+u_{\sigma_{n+1}\sigma_{n+1+s}}
u_{\sigma_k\sigma_{n+1+s}}
(\sigma_{n+1},\sigma_{n+1+s})(\sigma_k,\sigma_{n+1+s})
\cosh(Jv{\rho}_{{n+1-s}})-
$$
$$
-i\left [ u_{\sigma_{n+1}\sigma_{s'}}
u_{\sigma_k\sigma_s}
(\sigma_{n+1},\sigma_{s'})(\sigma_k,\sigma_s)
\sinh(Jv{\rho}_{{s}})- \right.
$$
  \begin{equation}
\left. \left. -u_{\sigma_{n+1}\sigma_{n+1-s}}
u_{\sigma_k\sigma_{n+1+s}}
(\sigma_{n+1},\sigma_{n+1-s})(\sigma_k,\sigma_{n+1+s})
 \sinh(Jv{\rho}_{{n+1-s}}) \right ] \right\}=0,
\label{o1-7}
\end{equation}
$$
u_{\sigma_{n+1}\sigma_{n+1}}u_{\sigma_{n+1+k}\sigma_{n+1}}
(\sigma_{n+1+k},\sigma_{n+1})+
$$
$$
+\sum^n_{s=1}
\left\{  u_{\sigma_{n+1}\sigma_s}u_{\sigma_{n+1+k}\sigma_s}
(\sigma_{n+1},\sigma_s)(\sigma_{n+1+k},\sigma_s)
 \cosh(Jv{\rho}_{{s}})+ \right.
$$
$$
+u_{\sigma_{n+1}\sigma_{n+1+s}}
u_{\sigma_{n+1+k}\sigma_{n+1+s}}
(\sigma_{n+1},\sigma_{n+1+s})(\sigma_{n+1+k},\sigma_{n+1+s})
 \cosh(Jv{\rho}_{{n+1-s}})-
$$
$$
-i \left [ 
u_{\sigma_{n+1}\sigma_{s'}}u_{\sigma_{n+1+k}\sigma_s}
(\sigma_{n+1+k},\sigma_s)(\sigma_{n+1},\sigma_{s'})
\sinh(Jv{\rho}_{{s}})- \right.
$$
  \begin{equation}
\left. \left.
-u_{\sigma_{n+1}\sigma_{n+1-s}}u_{\sigma_{n+1+k}\sigma_{n+1+s}}
(\sigma_{n+1+k},\sigma_{n+1+s})(\sigma_{n+1},\sigma_{n+1-s})
\sinh(Jv{\rho}_{{n+1-s}}) \right ] \right\}=0,
\label{o1-8}
\end{equation}
$$
u_{\sigma_{n+1+k}\sigma_{n+1}}
u_{\sigma_{n+1}\sigma_{n+1}}
(\sigma_{n+1+k},\sigma_{n+1})+
$$
$$
+\sum^n_{s=1}
\left\{ u_{\sigma_{n+1+k}\sigma_s}
u_{\sigma_{n+1}\sigma_s}
(\sigma_{n+1},\sigma_{s})(\sigma_{n+1+k},\sigma_{s})
\cosh(Jv{\rho}_{{s}})+ \right.
$$
$$
+u_{\sigma_{n+1+k}\sigma_{n+1+s}}
u_{\sigma_{n+1}\sigma_{n+1+s}}
(\sigma_{n+1+k},\sigma_{n+1+s})(\sigma_{n+1},\sigma_{n+1+s})
\cosh(Jv{\rho}_{{n+1-s}})+
$$
$$
 +i\left [
u_{\sigma_{n+1+k}\sigma_s}
u_{\sigma_{n+1}\sigma_{s'}}
(\sigma_{n+1+k},\sigma_s)(\sigma_{n+1},\sigma_{s'})
\sinh(Jv{\rho}_{{s}})- \right.
$$
\begin{equation}
\left. \left. -u_{\sigma_{n+1+k}\sigma_{n+1+s}}
u_{\sigma_{n+1}\sigma_{n+1-s}}
(\sigma_{n+1+k},\sigma_{n+1+s})
(\sigma_{n+1},\sigma_{n+1-s})
\sinh(Jv{\rho}_{{n+1-s}}) \right ] \right\}=0,
\label{o1-9}
\end{equation}
and additional relations
$ U^t(j;\sigma)\tilde{C}^{-1}_v(j)U(j;\sigma)=\tilde{C}^{-1}_v(j) $
 are equal
$$
u_{\sigma_{n+1}\sigma_k}
u_{\sigma_{n+1}\sigma_p}
(\sigma_{n+1},\sigma_k)(\sigma_{n+1},\sigma_p)
+\sum^n_{s=1}
\left\{ u_{\sigma_s\sigma_k}
u_{\sigma_s\sigma_p}
(\sigma_s,\sigma_k)(\sigma_s,\sigma_p)
\cosh(Jv{\rho}_{{s}})+ \right.
$$
$$
+u_{\sigma_{n+1+s}\sigma_k}
u_{\sigma_{n+1+s}\sigma_p}
(\sigma_{n+1+s},\sigma_k)(\sigma_{n+1+s},\sigma_p)
\cosh(Jv{\rho}_{{n+1-s}})+
$$
$$
 +i\left [ u_{\sigma_{s'}\sigma_k}
 u_{\sigma_s\sigma_p}
(\sigma_{s'},\sigma_k)(\sigma_s,\sigma_p)
\sinh(Jv{\rho}_{{s}})-  \right.
$$
$$
\left. \left.-u_{\sigma_{n+1-s}\sigma_k}
u_{\sigma_{n+1+s}\sigma_p}
(\sigma_{n+1-s},\sigma_k)(\sigma_{n+1+s},\sigma_p)
\sinh(Jv{\rho}_{{n+1-s}}) \right ]\right\}=
\delta_{kp}\cosh(Jv{\rho}_{{k}}),
$$
$$
u_{\sigma_{n+1}\sigma_{n+1+k}}
u_{\sigma_{n+1}\sigma_{n+1+p}}
(\sigma_{n+1},\sigma_{n+1+k})
(\sigma_{n+1},\sigma_{n+1+p})+
$$
$$
+\sum^n_{s=1}
\left\{ u_{\sigma_s\sigma_{n+1+k}}
u_{\sigma_s\sigma_{n+1+p}}
(\sigma_s,\sigma_{n+1+k})(\sigma_s,\sigma_{n+1+p})
\cosh(Jv{\rho}_{{s}})+  \right.
$$
$$
+u_{\sigma_{n+1+s}\sigma_{n+1+k}}
u_{\sigma_{n+1+s}\sigma_{n+1+p}}
(\sigma_{n+1+s},\sigma_{n+1+k})
(\sigma_{n+1+s},\sigma_{n+1+p})
\cosh(Jv{\rho}_{{n+1-s}})+
$$
$$
+ i\left [ u_{\sigma_{s'}\sigma_{n+1+k}}
u_{\sigma_s\sigma_{n+1+p}}
(\sigma_{s'},\sigma_{n+1+k})
(\sigma_s,\sigma_{n+1+p})
\sinh(Jv{\rho}_{{s}})- \right.
$$
$$
\left. \left.-u_{\sigma_{n+1-s}\sigma_{n+1+k}}
u_{\sigma_{n+1+s}\sigma_{n+1+p}}
(\sigma_{n+1-s},\sigma_{n+1+k})
(\sigma_{n+1+s},\sigma_{n+1+p})
\sinh(Jv{\rho}_{{n+1-s}}) \right ] \right\}=
$$
$$
=\delta_{kp}\cosh(Jv{\rho}_{{n+1-k}}),
$$
$$
u_{\sigma_{n+1}\sigma_k}
u_{\sigma_{n+1}\sigma_{n+1+p}}
(\sigma_{n+1},\sigma_k)
(\sigma_{n+1},\sigma_{n+1+p})+
$$
$$
+\sum^n_{s=1}
\left\{ u_{\sigma_s\sigma_k}
u_{\sigma_s\sigma_{n+1+p}}
(\sigma_s,\sigma_k)
(\sigma_s,\sigma_{n+1+p})
\cosh(Jv{\rho}_{{s}})+   \right.
$$
$$
+u_{\sigma_{n+1+s}\sigma_k}
u_{\sigma_{n+1+s}\sigma_{n+1+p}}
(\sigma_{n+1+s},\sigma_k)
(\sigma_{n+1+s},\sigma_{n+1+p})
\cosh(Jv{\rho}_{{n+1-s}})+
$$
$$
+i\left [u_{\sigma_{s'}\sigma_k}
u_{\sigma_s\sigma_{n+1+p}}
(\sigma_{s'},\sigma_k)
(\sigma_s,\sigma_{n+1+p})
\sinh(Jv{\rho}_{{s}})- \right.
$$
$$
\left. \left.-u_{\sigma_{n+1-s}\sigma_k}
u_{\sigma_{n+1+s}\sigma_{n+1+p}}
(\sigma_{n+1-s},\sigma_k)
(\sigma_{n+1+s},\sigma_{n+1+p})
\sinh(Jv{\rho}_{{n+1-s}}) \right ]\right\}=
$$
$$
=-i\delta_{n+1-k,p}\sinh(Jv{\rho}_{{k}}),
$$
$$
u_{\sigma_{n+1}\sigma_{n+1+k}}
u_{\sigma_{n+1}\sigma_p}
(\sigma_{n+1},\sigma_{n+1+k})
(\sigma_{n+1},\sigma_p)+
$$
$$
+\sum^n_{s=1}
\left\{ u_{\sigma_s\sigma_{n+1+k}}
u_{\sigma_s\sigma_p}
(\sigma_s,\sigma_{n+1+k})
(\sigma_s,\sigma_p)
\cosh(Jv{\rho}_{{s}})+ \right.
$$
$$
+u_{\sigma_{n+1+s}\sigma_{n+1+k}}
u_{\sigma_{n+1+s}\sigma_p}
(\sigma_{n+1+s},\sigma_{n+1+k})
(\sigma_{n+1+s},\sigma_p)
\cosh(Jv{\rho}_{{n+1-s}})+
$$
$$
+i\left [u_{\sigma_{s'}\sigma_{n+1+k}}
u_{\sigma_s\sigma_p}
(\sigma_{s'},\sigma_{n+1+k})
(\sigma_s,\sigma_p)
\sinh(Jv{\rho}_{{s}})- \right.
$$
$$
\left. \left.-u_{\sigma_{n+1-s}\sigma_{n+1+k}}
u_{\sigma_{n+1+s}\sigma_p}
(\sigma_{n+1-s},\sigma_{n+1+k})
(\sigma_{n+1+s},\sigma_p)
\sinh(Jv{\rho}_{{n+1-s}}) \right ]\right\} =
$$
$$
=i \delta_{n+1-k,p}\sinh(Jv{\rho}_{{p}}),
$$
$$
u^2_{\sigma_{n+1}\sigma_{n+1}}+ \sum^n_{k=1}
\left\{ u^2_{\sigma_k\sigma_{n+1}}
(\sigma_k,\sigma_{n+1})^2
\cosh(Jv{\rho}_{{k}})+ \right.
$$
$$
+u^2_{\sigma_{n+1+k}\sigma_{n+1}}
(\sigma_{n+1+k},\sigma_{n+1})^2
\cosh(Jv{\rho}_{{n+1-k}})+ 
$$
$$
+i\left [u_{\sigma_{n+1+k}\sigma_{n+1}}
u_{\sigma_{n+1-k}\sigma_{n+1}}
(\sigma_{n+1+k},\sigma_{n+1})
(\sigma_{n+1-k},\sigma_{n+1})
\sinh(Jv{\rho}_{{n+1-k}})- \right.
$$
$$
\left. \left. -u_{\sigma_k\sigma_{n+1}}
u_{\sigma_{k'}\sigma_{n+1}}
(\sigma_k,\sigma_{n+1})
(\sigma_{k'},\sigma_{n+1})
\sinh(Jv{\rho}_{{k}}) \right ] \right\}=1,
$$
$$
u_{\sigma_{n+1}\sigma_k}
u_{\sigma_{n+1}\sigma_{n+1}}
(\sigma_{n+1},\sigma_k)+
\sum^n_{p=1}
\left\{ u_{\sigma_p\sigma_k}
u_{\sigma_p\sigma_{n+1}}
(\sigma_p,\sigma_k)
(\sigma_p,\sigma_{n+1})
\cosh(Jv{\rho}_{{p}})+ \right.
$$
$$
+u_{\sigma_{n+1+p}\sigma_k}
u_{\sigma_{n+1+p}\sigma_{n+1}}
(\sigma_{n+1+p},\sigma_k)
(\sigma_{n+1+p},\sigma_{n+1})
\cosh(Jv{\rho}_{{n+1-p}})-
$$
$$
-i\left [u_{\sigma_p\sigma_k}
u_{\sigma_{p'}\sigma_{n+1}}
(\sigma_p,\sigma_k)
(\sigma_{p'},\sigma_{n+1})
\sinh(Jv{\rho}_{{p}})- \right.
$$
$$
\left. \left.-u_{\sigma_{n+1+p}\sigma_k}
u_{\sigma_{n+1-p}\sigma_{n+1}}
(\sigma_{n+1+p},\sigma_k)
(\sigma_{n+1-p},\sigma_{n+1})
\sinh(Jv{\rho}_{{n+1-p}}) \right ] \right\}=0,
$$
$$
u_{\sigma_{n+1}\sigma_{n+1}}
u_{\sigma_{n+1}\sigma_k}
(\sigma_{n+1},\sigma_k)+
$$
$$
+\sum^n_{p=1}
\left\{ u_{\sigma_p\sigma_{n+1}}
u_{\sigma_p\sigma_k}
(\sigma_p,\sigma_{n+1})
(\sigma_p,\sigma_k)
\cosh(Jv{\rho}_{{p}})+  \right.
$$
$$
+u_{\sigma_{n+1+p}\sigma_{n+1}}
u_{\sigma_{n+1+p}\sigma_k}
(\sigma_{n+1+p},\sigma_{n+1})
(\sigma_{n+1+p},\sigma_k)
\cosh(Jv{\rho}_{{n+1-p}})+
$$
$$
 +i\left [u_{\sigma_{p'}\sigma_{n+1}}
 u_{\sigma_p\sigma_k}
 (\sigma_{p'},\sigma_{n+1})
(\sigma_p,\sigma_k)
 \sinh(Jv{\rho}_{{p}})- \right.
$$
$$
\left. \left.-u_{\sigma_{n+1+p}\sigma_{n+1}}
 u_{\sigma_{n+1+p}\sigma_k}
 (\sigma_{n+1+p},\sigma_{n+1})
(\sigma_{n+1+p},\sigma_k)
 \sinh(Jv{\rho}_{{n+1-p}}) \right ] \right\}=0,
$$
$$
u_{\sigma_{n+1}\sigma_{n+1}}
u_{\sigma_{n+1}\sigma_{n+1+k}}
(\sigma_{n+1},\sigma_{n+1+k})+
$$
$$
+\sum^n_{p=1}
\left\{u_{\sigma_p\sigma_{n+1}}
u_{\sigma_p\sigma_{n+1+k}}
(\sigma_p,\sigma_{n+1})
(\sigma_p,\sigma_{n+1+k})
\cosh(Jv{\rho}_{{p}})+   \right.
$$
$$
+u_{\sigma_{n+1+p}\sigma_{n+1}}
u_{\sigma_{n+1+p}\sigma_{n+1+k}}
(\sigma_{n+1+p},\sigma_{n+1})
(\sigma_{n+1+p},\sigma_{n+1+k})
 \cosh(Jv{\rho}_{{n+1-p}})+
 $$
 $$
+i\left [u_{\sigma_{p'}\sigma_{n+1}}
u_{\sigma_p\sigma_{n+1+k}}
(\sigma_{p'},\sigma_{n+1})
(\sigma_p,\sigma_{n+1+k})
\sinh(Jv{\rho}_{{p}})- \right.
$$
$$
\left. \left.-u_{\sigma_{n+1-p}\sigma_{n+1}}
u_{\sigma_{n+1+p}\sigma_{n+1+k}}
(\sigma_{n+1-p},\sigma_{n+1})
(\sigma_{n+1+p},\sigma_{n+1+k})
\sinh(Jv{\rho}_{{n+1-p}}) \right ] \right\}=0,
$$
$$
u_{\sigma_{n+1}\sigma_{n+1+k}}
u_{\sigma_{n+1}\sigma_{n+1}}
(\sigma_{n+1},\sigma_{n+1+k})+
$$
$$
+\sum^n_{p=1}
\left\{ u_{\sigma_p\sigma_{n+1+k}}
u_{\sigma_p\sigma_{n+1}}
(\sigma_p,\sigma_{n+1+k})
(\sigma_p,\sigma_{n+1})
\cosh(Jv{\rho}_{{p}})+ \right.
$$
$$
+u_{\sigma_{n+1+p}\sigma_{n+1+k}}
u_{\sigma_{n+1+p}\sigma_{n+1}}
(\sigma_{n+1+p},\sigma_{n+1+k})
(\sigma_{n+1+p},\sigma_{n+1})
 \cosh(Jv{\rho}_{{n+1-p}})-
 $$
$$
- i\left [u_{\sigma_p\sigma_{n+1+k}}
u_{\sigma_{p'}\sigma_{n+1}}
(\sigma_p,\sigma_{n+1+k})
(\sigma_{p'},\sigma_{n+1})
\sinh(Jv{\rho}_{{p}})- \right.
 $$
 \begin{equation}
\left. \left. -u_{\sigma_{n+1+p}\sigma_{n+1+k}}
u_{\sigma_{n+1-p}\sigma_{n+1}}
(\sigma_{n+1+p},\sigma_{n+1+k})
(\sigma_{n+1-p},\sigma_{n+1})
\sinh(Jv{\rho}_{{n+1-p}}) \right ] \right\}=0,
\label{ort-1}
\end{equation}
where $ k,p=1, \ldots ,n. $  
$(v,j)$-orthogonality relations of  quantum group 
$ SO_v(N;j;\sigma),$ $ N=2n $
are given by above-mentioned formulae with the replacement
$n+1$ on $n.$


\end{document}